\documentclass{article}
\usepackage{arxiv}

\usepackage[utf8]{inputenc} % allow utf-8 input
\usepackage[T1]{fontenc}    % use 8-bit T1 fonts
\usepackage{url}            % simple URL typesetting
\usepackage{booktabs}       % professional-quality tables
\usepackage{nicefrac}       % compact symbols for 1/2, etc.
\usepackage{microtype}      % microtypography

\usepackage[english]{babel}
\usepackage{caption}
\usepackage{subcaption}
\usepackage{graphics}

\usepackage{amsfonts,amssymb,amsmath,amsthm}
\usepackage{breqn}
\usepackage{epsfig}
\usepackage{epstopdf}
\usepackage{amstext}
\usepackage{amsbsy}
\usepackage{amsopn}
\usepackage{float}
\usepackage{lipsum}
\usepackage{appendix}

\usepackage{hyperref}
\usepackage{cleveref}
\usepackage{setspace}

\newcommand\pd[3]{\frac{\partial^{#1} #2}{\partial #3^{#1}}}
\newcommand\status{false}

\newenvironment{rcases}
  {\left.\begin{aligned}}
  {\end{aligned}\right\rbrace}

\title{Combined homotopy and Galerkin stability analysis of Mathieu-like equations}

\author{
  Jeet Desai\\
  Universit\'e Pierre et Marie Curie\\
  Paris, France \\
  \texttt{jeet.desai@etu.upmc.fr } \\
   \And
 Amol Marathe \\
  Department of Mechanical Engineering\\
  BITS-Pilani\\
  Pilani, Rajasthan \\
  \texttt{amolmm@pilani.bits-pilani.ac.in } \\}

\begin{document}
\maketitle

\delimiterfactor=550 

\begin{abstract}
We propose homotopy analysis method in combination with Galerkin projections to obtain transition curves of Mathieu-like equations. While constructing homotopy, we think of convergence-control parameter as a function of embedding parameter and call it a convergence-control function. Homotopy analysis provides a relation between the parameters of Mathieu equation that also includes free parameters arising from convergence-control function. We generate extra nonlinear algebraic equations using Galerkin projections and solve numerically for arriving at transition curves. We demonstrate the usefulness of our method in case of three distinct versions of linear Mathieu equation by carefully choosing nonlinear and auxiliary  operators. Since homotopy analysis does not demand smallness of any of the parameters, our approach has distinct advantage over perturbation methods in determining transition curves covering a large region of the parameter space. The method is applicable to a wide variety of parametrically excited oscillators. 
\end{abstract}

% keywords can be removed
\keywords{ Mathieu\and Transition curves\and Homotopy analysis\and Galerkin projections}

\section{Introduction}

The well-known Mathieu equation models diverse physical phenomena including the motion of a simple pendulum whose support is excited periodically in vertical direction as well as the motion of ions within Paul trap mass spectrometer. It is a linear, second order ordinary differential equation with time periodic coefficient and is given by
\begin{equation*}
    \ddot{x} + \Big(\delta  + \epsilon \cos(t)\Big) x = 0. 
\end{equation*}
A chosen point in the $\delta$-$\epsilon$ plane is stable or unstable depending on whether solutions remain bounded with time or grow unboundedly. Stable regions in the parameter plane can be generated using numerical calculation of Floquet multipliers at grid points (see Appendix \ref{floquet}). Transition curves are nothing but the boundaries of the stable regions, which are not just typical of linear Mathieu equation, but also occur in case of several Mathieu-like oscillators, for example, damped Mathieu equation, quasi-periodic Mathieu, fractional-order Mathieu and Mathieu with excitation of non-circular type. Finding transition curves analytically is computationally advantageous over generating stable and unstable regions. Producing these curves has been attempted using the method of harmonic balance via Hill's determinants as one of the two linearly independent solutions to Mathieu equation is periodic along the transition curves. However, the method faces problems computationally due to the infinite size of the sets of linear homogeneous algebraic equations involved. Approximate transition curves have also been obtained using singular perturbation methods with $\epsilon$ as the perturbation parameter. Usually, these curves are restricted to a small region of $\delta$-$\epsilon$ plane as the validity of such methods is restricted to oscillators involving a small parameter. In this paper, we demonstrate the homotopy analysis method augmented by Galerkin projections seeking transition curves of three versions of Mathieu equation; classical Mathieu, linearly damped Mathieu and Mathieu with impulsive excitation. The approach developed here has the potential of obtaining transition curves covering a large region of the parameter space for  numerous Mathieu-like oscillators, including scalable Mathieu equations for which Floquet theory is not applicable. This technique can also be applied to get the steady-state response of externally forced strongly nonlinear oscillators.

Transition curves have been investigated by many researchers for different Mathieu-like equations over the past few decades. The connection between stability diagrams for linear Mathieu equation and second order delay differential equation was used for arriving at the same for delayed Mathieu equation \cite{Insperger1989}. Approximate steady-state response and transition curves of fractional Mathieu-Duffing equation were studied \cite{leung} using harmonic balance and polynomial homotopy continuation, a technique to solve a system of nonlinear algebraic equations. Stability diagrams and transition curves were studied numerically for a class of Mathieu equations for which principle of superposition is invalid; but the scalability of solutions holds \cite{Marathe1643}. Quasi-periodic Mathieu equation was studied for transition curves \cite{quasi} using regular perturbation and harmonic balance-Hill’s determinants. The parameter plane studied is not the usual $\delta$-$\epsilon$ plane, but $\delta$-$\omega,$ $\omega$ being an irrational frequency of a component of parametric excitation. Stability charts were generated for a pendulum with parametric excitation of the elliptic type using Floquet theory and harmonic balance \cite{Sah3995}. Analytical estimates of Floquet exponents were attempted for arbitrarily chosen parameter values and used to find transition curves \cite{gizem}. Field inhomogeneties alter stability boundaries in nonlinear Paul trap mass spectrometers. Changes in transition curves due to hexapole and octopole superposition were studied using harmonic balance method and continued fractions approach \cite{Sevugarajan2002181}. Early and delayed ejection of ions in the mass selective ejection experiments due to field inhomogeneties were modeled as dynamics in the neighborhood of the stability boundary of linear Mathieu equation and studied using the method of multiple scales \cite{Rajanbabu2007170}.

We apply homotopy analysis method (HAM) and Galerkin projections to get transition curves of Mathieu equation covering a large region of the parameter plane. HAM is an approximate analytical method that works well even in the absence of a small parameter and Galerkin projections is a heuristic technique used to approximate solutions of initial and boundary value problems. HAM is based on the concept of homotopy, a continuous deformation that takes a chosen initial function to the desired final function. This deformation is parametrized by embedding parameter $p,$ varying continuously from $0$ to $1.$ In the present context, a periodic solution to Mathieu equation is the desired function and is to be obtained by continuously deforming a solution of a simple harmonic oscillator, say, $\cos(t)$ or $\sin(t)$ or a combination of the two. Accordingly, we construct the homotopy. In case of Mathieu equation with damping, we modify the homotopy with an additional term. Such flexibility with HAM is of paramount importance. Though homotopy analysis is not a true perturbation method, the solution here is to be obtained in the spirit of a typical singular perturbation method: HAM breaks down the original nonlinear oscillator into a sequence of externally forced linear oscillators, except at zeroth order. Also while obtaining a bounded solution, HAM usually involves the important step of the removal of unbounded terms from solutions at every order. To adjust the convergence region of the series solution at $p=1,$ the convergence-control parameter $h$ is introduced into the homotopy that is usually considered to be a constant. In our application, we think of this parameter as an unknown function of the embedding parameter $p,$ calling it a convergence-control function. Different derivatives of this function evaluated at $p=0$ appear in the series solution acting as free parameters which when found suitably improves the convergence substantially. The removal of secular terms upto order $N+1$ yields one equation in $\delta$ and $\epsilon$  relating period of the solution of the Mathieu equation to the period of the solution of a simple harmonic oscillator. This equation also involves $N+1$ free parameters mentioned above. Choosing $N+1$ appropriate weighting functions, Galerkin projections yield equal number of equations. Numerically solving the resulting system of simultaneous algebraic equations yields a curve in $N+3$-dimensional space (free parameters, $\delta$ and $\epsilon$) whose projection in $\delta$-$\epsilon$ plane is nothing but the approximate transition curve.

The remainder of this paper is structured as follows. In section \ref{sec:2}, we describe the method, as applied to linear Mathieu equation and discuss a few associated numerical issues along with the results obtained. We modify the method in certain ways in section \ref{sec:3} to make it applicable to the Mathieu equation with damping. We amend the method again for Mathieu equation with impulsive parametric excitation to ease the computational effort in section \ref{sec:4}. The final section draws our work to a close, while pointing at directions for further work.

\section{Linear Mathieu equation}
\label{sec:2}

Consider a linear Mathieu equation,
\begin{equation} \label{ode_lin_math}
	\ddot{x} + \Big(\delta  + \epsilon \cos(t)\Big) x = 0. 
\end{equation}
We wish to obtain a periodic solution $x(t)$ to Eqn.(\ref{ode_lin_math}) using homotopy analysis method. As usual with HAM, we introduce an embedding parameter $p \in [0,1]$ and function $\tilde{x}(t;p)$ for constructing homotopy as given by
\begin{equation}\label{PDE_lin_math_unscaled}
	\mathcal{H} \equiv (1-p)\mathcal{L}\Big(\tilde{x}(t;p) \Big) - h(p) \mathcal{N}\Big(\tilde{x}(t;p) \Big) = \ 0. 
\end{equation}
Here, $\mathcal{L}$  and $\mathcal{N}$ are linear and nonlinear operators respectively. As $p$ varies continuously from $0$ to $1,$ the solution of homotopy $\mathcal{H}$ varies continuously from the solution of $\mathcal{L}\Big(\tilde{x}(t;0)\Big)=0$ to the solution of  $\mathcal{N}\Big(\tilde{x}(t;1)\Big)=0.$ We set $\mathcal{L}\Big(\tilde{x}(t;0)\Big)\equiv \ \ddot{x}+x,$ a simple harmonic oscillator and $\mathcal{N}\Big(\tilde{x}(t;1)\Big)$ to Eqn.(\ref{ode_lin_math}). In the homotopy framework, the desired periodic solution is $\tilde{x}(t;1).$

While constructing the homotopy $\mathcal{H},$ we differ from the usual HAM in a small way. Coefficient of nonlinear operator $\mathcal{N}$ in Eqn.(\ref{PDE_lin_math_unscaled}) is usually the product of convergence-control and embedding parameters, i.e., $hp$ \cite{Wen2007427}, which we replace by a function, say convergence-control function $h(p)$ subject to conditions $h(0)=0$ and $h(1) \neq 0$. This alteration ensures much faster convergence of the solution of Eqn.(\ref{PDE_lin_math_unscaled}). However, derivatives of $h(p)$ at $p=0$ appear as unknowns in solutions at different orders.

With $\omega(1)$ as the frequency of periodic solution of Mathieu equation, we now scale time $t$ to $\tau=\omega(p)t$ and introduce time-stretching function $\lambda(p)=\tfrac{1}{\omega(p)}.$ We also choose $\lambda(0)=1.$  Depending on whether we are looking for periodic solution with period $2\pi$ or $4\pi,$ we set $\lambda(1)=1$ or $2$ respectively and proceed with $\lambda(1)=2.$ Homotopy (Eqn.(\ref{PDE_lin_math_unscaled})) with respect to the scaled time is
\begin{equation}\label{PDE_lin_math_scaled}
\mathcal{H} \equiv (1-p) \Big( \tilde{x}_{\tau\tau} + \lambda(p)^{2}\tilde{x} \Big)
 - h(p)\bigg( \tilde{x}_{\tau\tau} + \lambda(p)^{2}\Big(\delta + \epsilon \cos\big(\lambda(p)\tau\big)\Big) \tilde{x} \bigg)=0,
\end{equation}
where subscript $\tau$ denotes differentiation w.r.t. $\tau.$ We are free to alter Eqn.\eqref{PDE_lin_math_scaled} provided it matches with $\tau$-scaled version of Eqn.\eqref{ode_lin_math} at $p=1.$ We introduce the new homotopy
\begin{equation}\label{PDE_lin_math_scaled_lambda}
\mathcal{H} \equiv (1-p) \Big( \tilde{x}_{\tau\tau} + \lambda(p)^{2}\tilde{x} \Big)
 - h(p)\bigg( \tilde{x}_{\tau\tau} + \lambda(p)^{2}\Big(\delta + \epsilon \cos\big(\lambda(1)\tau\big)\Big) \tilde{x} \bigg)=0.
\end{equation}
Noting that Eqn.\eqref{PDE_lin_math_scaled} and Eqn.\eqref{PDE_lin_math_scaled_lambda} coincide at $p=1,$ we will find it more convenient to work with Eqn.\eqref{PDE_lin_math_scaled_lambda}. Taylor-expanding $\tilde{x}(\tau;p)$ about $p=0,$ we have
\begin{align}
	\tilde{x}(\tau;p) &= x^{\mbox{\tiny [0]}}(\tau) + \sum_{n=1}^{\infty}\frac{1}{n!} x^{\mbox{\tiny [{\em n}]}}(\tau)p^{n},\label{solution_taylor}\\
	\mbox{where}  \quad  x^{\mbox{\tiny [0]}}(\tau) &= \tilde{x}(\tau;0)  \quad \mbox{and} \quad  x^{\mbox{\tiny [{\em n}]}}(\tau) = \left.\pd{n}{\tilde{x}}{p} \right|_{p=0}. \nonumber
\end{align}
Differentiating Eqn.(\ref{solution_taylor}) w.r.t. $\tau$, we get
\begin{align}
   \tilde{x}_{\tau} = \dot{x}^{\mbox{\tiny [0]}}(\tau) + \sum_{n=1}^{\infty}\frac{1}{n!} \dot{x}^{\mbox{\tiny [{\em n}]}}(\tau)p^{n}.  \label{solution_taylor_xdot}
\end{align}
It is quite clear from the context that overdot in the above equation denotes differentiation w.r.t. $\tau$ and not w.r.t. $t.$  
We Taylor-expand time-stretching function and convergence-control function as well to get 

\begin{align}
	\lambda(p) = 1 + \sum_{n=1}^{\infty}\frac{1}{n!} \lambda^{\mbox{\tiny [{\em n}]}}p^{n} \quad &\mbox{where} \quad \lambda^{\mbox{\tiny [{\em n}]}}=  \left.\frac{d^{n}\lambda(p)}{dp^{n}} \right|_{p=0}\label{solution_taylor_lambda}\mbox{, and}\\
	h(p) =  \sum_{n=1}^{\infty}\frac{1}{n!} h^{\mbox{\tiny [{\em n}]}}p^{n} \quad &\mbox{where} \quad h^{\mbox{\tiny [{\em n}]}} =  \left.\frac{d^{n}h(p)}{dp^{n}} \right|_{p=0}.\label{solution_taylor_h}
\end{align}
We neither need to specify $h(p)$ as a function of $p$ nor need to find it. However, its derivatives at $p=0$ will appear as unknowns in the subsequent analysis and play a significant role in accelerating the convergence of the solution $\tilde{x}(t;1).$ We fix the number of terms in the Taylor expansion of $\tilde{x}(\tau;p)$ to a small positive integer, say $N.$ Typically, we will use $N=3.$

We begin by obtaining $x^{\mbox{\tiny [0]}}(\tau).$ Substituting $p=0$ in Eqn.(\ref{PDE_lin_math_scaled}), we get the zeroth-order deformation equation
\begin{equation}\label{zeroth_order_eqn_impulse}
	\mathcal{L}\big(x^{\mbox{\tiny[0]}}\big) = \ddot{x}^{\mbox{\tiny [0]}}+x^{\mbox{\tiny[0]}} = 0.
\end{equation}
For a linear Mathieu equation, the unstable regions emanate at $\tfrac{k^2}{4}\;(k=1,2,\cdots)$ from $\delta$-axis in the $\delta$-$\epsilon$ plane \cite{RRand}. Along the boundary of the alternate unstable regions, there exists periodic solutions with period $2\pi$ and $4\pi.$ The stability boundary of a given unstable region consists of the left branch and the right branch, also called transition curves. The periodic solutions along one branch can be obtained with the initial condition $(1,0),$ and along the other with $(0,1)$. In order to get a branch with the  initial condition $(1,0)$ we consider Eqns.(\ref{solution_taylor}) and (\ref{solution_taylor_xdot}), and set $\tilde{x}(0;1)=1$ and $\tilde{x}_{\tau}(0;1)=0.$ Accordingly, we choose the initial condition for zeroth-order deformation equation as $x^{\mbox{\tiny [0]}}(0)=1$ and $\dot{x}^{\mbox{\tiny [0]}}(0)=0$ to obtain
\begin{equation}\label{zeroth_order_solution} 
	x^{\mbox{\tiny [0]}}(\tau) = \cos(\tau).
\end{equation}
To get the first-order deformation equation, we differentiate Eqn.(\ref{PDE_lin_math_scaled_lambda}) w.r.t. $p$ once. Then substituting $p=0,$ and using Eqns.\eqref{solution_taylor},\eqref{solution_taylor_lambda},\eqref{solution_taylor_h} and \eqref{zeroth_order_solution}, we get
\begin{equation}\label{first_order}
	\ddot{x}^{\mbox{\tiny [1]}} + x^{\mbox{\tiny [1]}} = \bigg(-2 \lambda^{\mbox{\tiny [1]}} +\left(- 1 + \delta + \frac{\epsilon}{2} \right) h^{\mbox{\tiny [1]}} \bigg) \cos(\tau)+\frac{\epsilon h^{\mbox{\tiny [1]}}}{2}  \cos(3 \tau).
\end{equation}
Equation (\ref{first_order}) is an externally forced simple harmonic oscillator. We now remove the coefficient of resonant forcing, i.e., the secular term and thus determine
\[
	\lambda^{\mbox{\tiny [1]}} = \frac{h^{\mbox{\tiny [1]}}}{2} \left(- 1 + \delta + \frac{\epsilon}{2} \right).
\]
$h^{\mbox{\tiny [1]}}$ is yet undetermined and is playing a role of the free parameter in choosing curves in the parameter plane. Our choice of the initial condition for zeroth-order deformation equation fixes the initial conditions for the higher order deformation equations as
\begin{equation}\label{initial_conditions_nth_order_lin}
	x^{\mbox{\tiny[{\em n}]}}(0)=0, \quad \dot{x}^{\mbox{\tiny[{\em n}]}}(0)=0, \quad n\geq1. 
\end{equation}
Using Eqn.(\ref{initial_conditions_nth_order_lin}), we solve the first-order deformation equation to get 
\[ 
	x^{\mbox{\tiny [1]}}(\tau) = \frac{\epsilon h^{\mbox{\tiny [1]}}}{16}\Big( \cos \left(\tau\right)-\cos ( 3 \tau) \Big). 
\]
Differentiating Eqn.(\ref{PDE_lin_math_scaled_lambda}) w.r.t. $p$ twice and evaluating the same at $p=0,$  we get the second-order deformation equation. By substituting the solutions of zeroth- and first-order deformation equations in this, we get
\begin{multline}\label{second_order_eqn}
	\ddot{x}^{\mbox{\tiny[2]}} + x^{\mbox{\tiny[2]}}  = \Bigg(-2\lambda^{\mbox{\tiny [2]}} +  \Big(2h^{\mbox{\tiny[1]}}+h^{\mbox{\tiny[2]}}\Big)\bigg(-1 + \delta  + \frac{\epsilon}{2}\bigg) + \frac{{h^{\mbox{\tiny[1]}}}^2}{2}\bigg(-1-2\delta-\epsilon + 3{\delta}^{2}+3\epsilon\delta+\frac{5 {\epsilon}^{2}}{8}\bigg)\Bigg)\cos(\tau) \\+ \epsilon\bigg(h^{\mbox{\tiny [1]}} + \frac{h^{\mbox{\tiny [2]}}}{2} + \delta {h^{\mbox{\tiny [1]}}}^{2} + \frac{5\epsilon {h^{\mbox{\tiny [1]}}}^{2}}{8} \bigg) \cos(3\tau)-\frac{{\epsilon}^{2}{h^{\mbox{\tiny [1]}}}^{2}}{16} \cos \left( 5 \tau \right).
\end{multline}
Elimination of the secular term yields
\begin{equation}
	\lambda^{\mbox{\tiny[2]}} = \bigg(h^{\mbox{\tiny[1]}}+\frac{h^{\mbox{\tiny[2]}}}{2}\bigg)\bigg(- 1 + \delta + \frac{\epsilon}{2}\bigg) + \frac{{h^{\mbox{\tiny[1]}}}^2}{4}\bigg( - 1  - 2\delta - \epsilon + 3{\delta}^{2}  + 3\epsilon\delta + \frac{5 {\epsilon}^{2}}{8}\bigg)  .
\end{equation}
Solving Eqn.(\ref{second_order_eqn}) after substitution of $\lambda^{\mbox{\tiny[2]}}$ from the above, we obtain
\begin{multline*}
	x^{\mbox{\tiny[2]}}(\tau) = \frac{\epsilon}{8}\bigg(\Big(h^{\mbox{\tiny[1]}}+ \frac{h^{\mbox{\tiny[2]}}}{2} + {\frac {29 {\epsilon}{h^{\mbox{\tiny[1]}}}^{2}}{48}}+ \delta {h^{\mbox{\tiny[1]}}}^{2}\Big) \cos \left( \tau \right) \\ -\Big(h^{\mbox{\tiny[1]}} + \frac{h^{\mbox{\tiny[2]}}}{2} +{\frac {5 {\epsilon}{h^{\mbox{\tiny[1]}}}^{2}}{8}} + \delta{h^{\mbox{\tiny[1]}}}^{2}\Big) \cos \left( 3 \tau \right) + {\frac {{\epsilon}{h^{\mbox{\tiny[1]}}}^{2}}{48}} \cos \left( 5 \tau \right) \bigg).
\end{multline*}
In a similar fashion, we proceed till $N$-th order deformation equation by removing secular terms to yield expressions for $\lambda^{\mbox{\tiny[1]}}$ to $\lambda^{\mbox{\tiny[{\em N}+1]}}.$ We notice that the expressions for $\lambda^{\mbox{\tiny[{\em n}]}}$ and $x^{\mbox{\tiny[{\em n}]}}(\tau)$ are in terms of free parameters $h^{\mbox{\tiny[{\em k}]}}, \ k=1,\cdots,n$ as well as $\delta$ and $\epsilon.$ Considering first $N+2$ terms of Eqn.(\ref{solution_taylor_lambda}) and enforcing it to $\lambda(1),$ we get a relation
\begin{equation}
	\lambda(1) = 2 = 1 + \sum_{n=1}^{N+1}\frac{1}{n!} \lambda^{\mbox{\tiny [{\em n}]}}.\label{taylor_approx1}
\end{equation}
By substituting the expressions for $\lambda^{\mbox{\tiny \tiny[1]}}$ to $\lambda^{\mbox{\tiny[{\em N}+1]}},$ we get an approximate expression for $\lambda(1).$ For $N=3,$ the corresponding expression is given in Appendix \ref{lin_math_eqns}. We approximate the periodic solution $\tilde{x}(\tau;1)$ by considering first $N+1$ terms of Eqn.(\ref{solution_taylor}) after substituting $p=1$ and expressions for $x^{\mbox{\tiny[{\em n}]}}(\tau).$ By scaling back to time $t$ from $\tau=\tfrac{t}{\lambda(1)}=\tfrac{t}{2}$ and denoting the approximate solution by $x_N(t),$ we get
\begin{equation}\label{taylor_approx4} 
	x_N(t)=x^{\mbox{\tiny[0]}}(t) + \sum_{n=1}^{N}\frac{1}{n!} x^{\mbox{\tiny[{\em n}]}}(t).
\end{equation}

Thus homotopy analysis furnishes us with one relation (Eqn.(\ref{taylor_approx1})) between $\delta$ and $\epsilon,$ but also involving  $h^{\mbox{\tiny[1]}},\cdots,h^{\mbox{\tiny[{\em N}+1]}}$ as unknowns. Therefore, we need more equations involving these unknowns solving which we hope to get transition curves in the parameter plane. We use Petrov-Galerkin method of weighted residual for obtaining these equations. Using $x_N(t)$ from Eqn.(\ref{taylor_approx4}), we get the residual
\begin{equation*}
 \mathcal{R}_N(t) \equiv  \ddot{x}_N(t) + \Big(\delta  + \epsilon \cos(t)\Big) x_N(t). 
 \end{equation*}
We choose the number of weighting functions as being equal to the number of equations desired. As per Petrov-Galerkin scheme, taking the weighting functions same as the basis functions of ${x}_N(t),$  i.e., $\cos(\tfrac{t}{2}),\;\cos(\tfrac{3t}{2}),\ \cos(\tfrac{5t}{2}), \cdots,$ we create $N+1$ equations
\begin{equation}\label{galerkin_lin_math}
	\int_{0}^{4\pi} \cos\Big((k+\tfrac{1}{2})t\Big)\mathcal{R}_N(t) dt = 0, \quad k=0,\cdots,N.
\end{equation}
Solving the above $N+1$ algebraic equations and Eqn.(\ref{taylor_approx1}) simultaneously for unknowns $\delta,\ \epsilon,\ h^{\mbox{\tiny[1]}},\cdots,h^{\mbox{\tiny[{\em N}+1]}}$ numerically, we obtain transition curves. 

\begin{figure}[!h]
\makebox[\textwidth][c]
{\includegraphics[draft=\status , width=15cm,height=14cm]{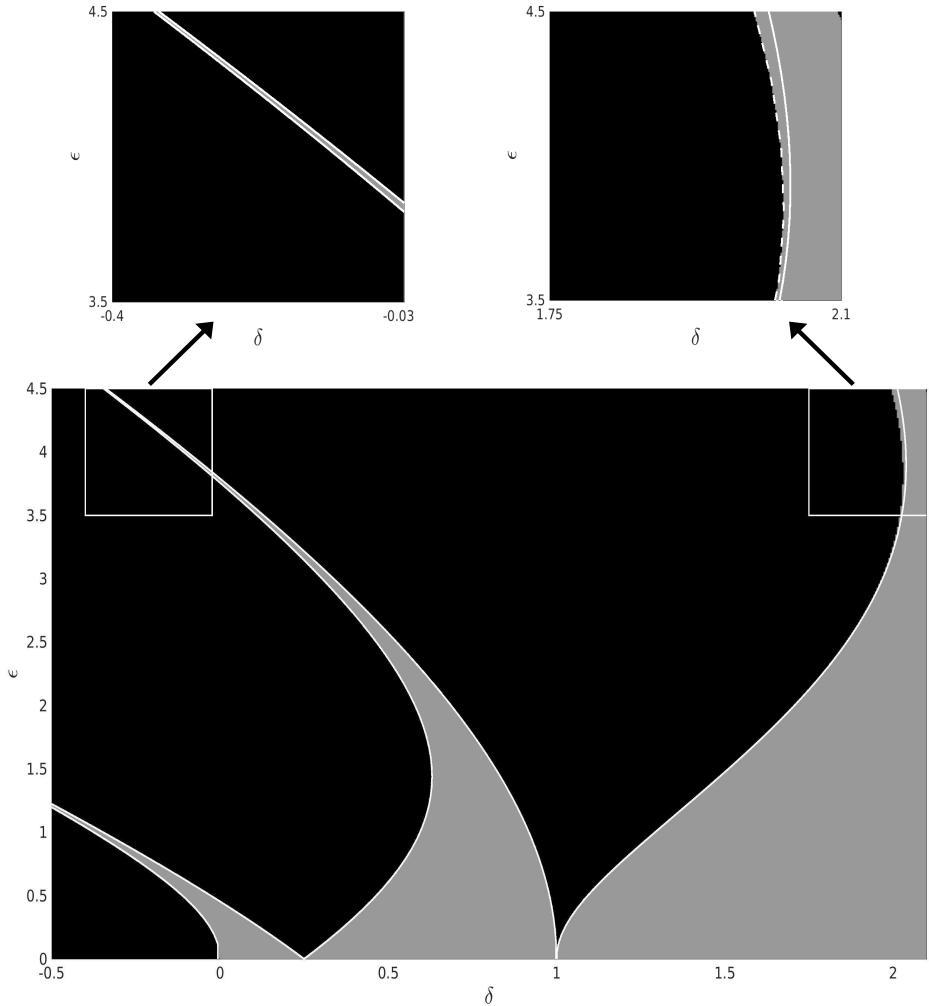}}
\caption{Stable regions (grey), unstable regions (black) and transition curves (white lines, homotopy$+$Galerkin)  in the parameter plane.}
\label{Linear_Mathieu_Stab}
\end{figure}

Figure (\ref{Linear_Mathieu_Stab}) shows the unstable regions in the parameter plane with $\delta$ varying from $-0.5$ to $2.1$ and $\epsilon$ from $0$ to $4.5.$ The region is divided into $780 \times 1350$ grid points. The Floquet multipliers are calculated numerically (Appendix \ref{floquet}) at each grid point and accordingly the stability is determined. $4\pi$ periodic solutions exist on the boundary of the unstable region emanating at $\delta=0.25.$ To get the left branch of the corresponding transition curve, we solve a system of five algebraic equations ($N=3,$ Appendix \ref{lin_math_eqns}) using fsolve command of a symbolic algebra package Maple. For finding a point on the transition curve, we fix $\epsilon$ and solve the equations for the remaining unknowns, including $\delta.$ We may get multiple $\delta$ values corresponding to different transition curves; but these can be eliminated by searching $\delta$ only within a pre-determined interval. The entire procedure is worked out with the  initial condition $(0,1)$ to get the right branch. Also plotted are the transition curves (white lines in Fig.(\ref{Linear_Mathieu_Stab})) originating from $\delta=0$ and $1$ corresponding to $2\pi$ periodic solutions ($\lambda(1)=1$).

\iffalse{
\begin{figure}[!h]
\makebox[\textwidth][c]
{\includegraphics[draft=\status , width=15cm,height=9.5cm]{Linear_Mathieu.eps}}
\caption{Stable regions (grey), unstable regions (black) and transition curves (white lines) in the parameter plane.}
\label{Linear_Mathieu_Stab}
\end{figure}

\begin{figure}[!h]
  \begin{subfigure}[b]{0.475\textwidth}
  \centering
    \includegraphics[draft=\status ,height=5cm,width=5cm]{Linear_Mathieu_zoom_2.eps}
    \caption{High accuracy with larger $\epsilon$}
    \label{fig:f1}
  \end{subfigure}
%  \hfill
  \begin{subfigure}[b]{0.475\textwidth}
\centering
    \includegraphics[draft=\status ,height=5cm,width=5cm]{Linear_Mathieu_zoom_1.eps}
    \caption{$N=3$ (conti. line) and $N=5$ (dashed). }
    \label{fig:f2}
  \end{subfigure}
  \caption{Zoomed portion of Fig \ref{Linear_Mathieu_Stab}. Grey, stable; black, unstable.}
  \label{zoom}
\end{figure}
\fi

\begin{figure}[!h]
\makebox[\textwidth][c]
{\includegraphics[draft=\status , width=18cm,height=11cm]{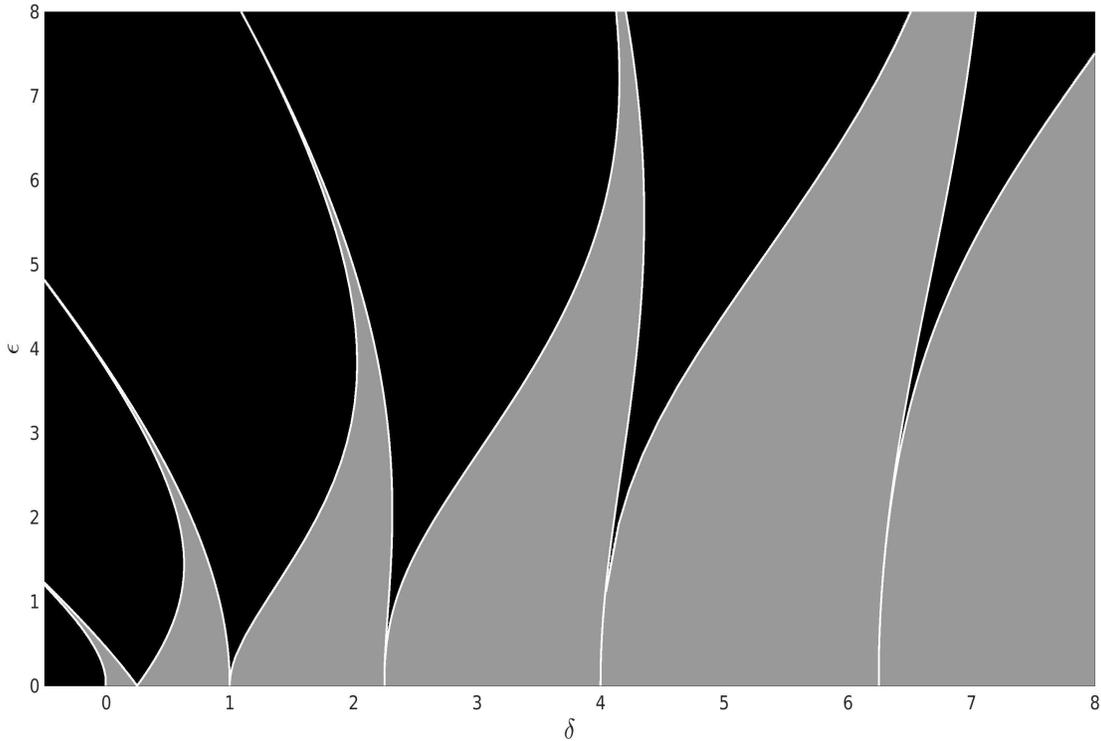}}
\caption{Stable regions (grey), unstable regions (black) and transition curves (white lines, homotopy$+$Galerkin) in the parameter plane.}
\label{Linear_Mathieu_Stab_Larger}
\end{figure}

The right branch of instability tongue originating from $\delta=0.25$ and the left branch of instability tongue originating from $\delta=1$ come closer in a small region highlighted using the left top white box (Fig.(\ref{Linear_Mathieu_Stab})). After zooming in, we observe that these branches as obtained for $N=3$ case are adequately separated. However, after zooming in the right top white box, it is clear that the transition curve with $N=3$ (continuous white line) does not capture the stability boundary with sufficient accuracy. We are able to capture the boundary to a  satisfactory extent with $N=5$ (dashed white line). We believe that transition curves corresponding to unstable regions emanating from $\delta=\tfrac{k^2}{4}$ for large $k,$ can be obtained for a fixed choice of $\lambda(0)=1.$ However, sufficiently high accuracy of the transition curves for such $\delta$ values can be maintained only with higher values of $N.$ This increases computational burden considerably. Based on our numerical experiments, we propose that the same level of accuracy can be maintained with smaller $N$ values, provided we choose larger value of $\lambda(0)=2,3$ or higher. Figure(\ref{Linear_Mathieu_Stab_Larger}) shows the transition curves obtained using $N=5$ extended over the larger region of the parameter plane depicting the usefulness of perturbation-like method.

\section{Linearly damped Mathieu equation} 
\label{sec:3}

Now we consider the Mathieu equation with linear damping,
\begin{equation} \label{ode_lin_damp_math}
 \ddot{x} + \Big(\delta  + \epsilon \cos(t)\Big) x + c\dot{x}= 0 .
\end{equation}
We choose the linear operator to be a simple harmonic oscillator as before and the nonlinear operator to be damped linear Mathieu equation to construct homotopy $\mathcal{H}$ as
\begin{equation*} 
\mathcal{H} \equiv (1-p)\Big(\tilde{x}_{tt} + \tilde{x}\Big) - h(p) \bigg( \tilde{x}_{tt} + \Big(\delta +\epsilon \cos(t)\Big) \tilde{x}+ c\tilde{x}_{t} \bigg) =0.
\end{equation*}
Scaling time $t$ to $\tau$ and replacing the term $\cos\Big(\lambda(p)\tau\Big)$ by $\cos\Big(\lambda(1)\tau\Big)$ as was done in the previous section, we get
\begin{equation}\label{PDE_lin_damp_scaled_lambda1}
\mathcal{H} \equiv  (1-p)\bigg( \tilde{x}_{\tau\tau} + \lambda(p)^2 \tilde{x}\bigg)  \\
-h(p) \bigg( \tilde{x}_{\tau\tau} + \lambda(p)^2\Big(\delta +\epsilon \cos\big(\lambda(1)\tau\big)\Big) \tilde{x} + \lambda(p)c\tilde{x}_{\tau}\bigg) =0.
\end{equation}
We choose again $\lambda(0)=1.$ This time we set $\lambda(1)=1$ to look for a periodic solution to Eqn.(\ref{ode_lin_damp_math}) with period $2\pi.$ In a routine manner, by substituting $p=0$ in the homotopy $\mathcal{H},$ we obtain the zeroth-order deformation equation
\begin{equation*}
 \mathcal{L}(x^{\mbox{\tiny [0]}}) = \ddot{x}^{\mbox{\tiny [0]}} + x^{\mbox{\tiny [0]}}=0.
\end{equation*}
For the Mathieu equation without damping, we know beforehand that   the periodic solutions exist along the transition curves with the  initial conditions $(1,0)$ or $(0,1),$ i.e., the initial phase zero or $\tfrac{\pi}{2}.$ However for Eqn.(\ref{ode_lin_damp_math}), we do not know at which phase to start with for arriving at the periodic solution for a given point on the transition curve. Also, since Eqn.(\ref{ode_lin_damp_math}) is scalable, the initial amplitude can be chosen arbitrarily to  get the required solution. Hence, we choose the initial condition as $(1,\zeta)$ with unknown $\zeta$ expressed as
\begin{equation*}
\zeta = \zeta^{\mbox{\tiny [0]}} + \sum_{n=1}^{\infty}\frac{1}{n!} \zeta^{\mbox{\tiny [{\em n}]}},
\end{equation*}
involving unknown $\zeta^{\mbox{\tiny [{\em n}]}}$ at $n$-th order. In an application of HAM to van der Pol oscillator, the  authors while choosing a convenient phase, expanded the unknown amplitude in a similar manner as above \cite{Chen20091816}.
Following the procedure from the previous section, we fix the initial conditions for $n$-th order deformation equation as 
\begin{align}\label{initial_conditions_nth_order_damp}
x^{\mbox{\tiny [{\em n}]}}(0)=0 \quad \mbox{and} \quad 
   \dot{x}^{\mbox{\tiny [{\em n}]}}(0)=\zeta^{\mbox{\tiny [{\em n}]}},  \quad n\geq 1.
\end{align}
For zeroth-order deformation equation, we set $x^{\mbox{\tiny [0]}}(0)=1$ and $\dot{x}^{\mbox{\tiny [0]}}(0)=\zeta^{\mbox{\tiny [0]}}$ to get
\begin{equation}\label{zeroth_sol_damp}
x^{\mbox{\tiny [0]}}(\tau) = \cos(\tau) + \zeta^{\mbox{\tiny [0]}}\sin(\tau).
\end{equation}
Substituting Eqn.(\ref{zeroth_sol_damp}) in first-order deformation equation, we get
\begin{multline*}
 \ddot{x}^{\mbox{\tiny [1]}} + x^{\mbox{\tiny [1]}}=
\Big(- 2\lambda^{\mbox{\tiny [1]}}+ \left(-1+\delta + c\zeta^{\mbox{\tiny[0]}}\right) h^{\mbox{\tiny[1]}}\Big) \cos(\tau) +\frac{\epsilon h^{\mbox{\tiny[1]}}}{2}\cos(2\tau)\\ - \Big( 2\zeta^{\mbox{\tiny[0]}} \lambda^{\mbox{\tiny[1]}}+ \left( c+
\zeta^{\mbox{\tiny[0]}}-\delta  \zeta^{\mbox{\tiny[0]}}\right) h^{\mbox{\tiny[1]}} \Big) \sin(\tau)+\frac{\epsilon \zeta^{\mbox{\tiny[0]}}h^{\mbox{\tiny[1]}}}{2}\sin(2\tau)+\frac{\epsilon h^{\mbox{\tiny[1]}}}{2}.
\end{multline*}
Eliminating the secular terms in the above equation, i.e., equating  coefficients of $\cos(\tau)$ and $\sin(\tau)$ to zero, we get two equations in unknowns $\lambda^{\mbox{\tiny[1]}}$ and $\zeta^{\mbox{\tiny[0]}}.$ By rearranging the two equations, we find that 
\[ {\zeta^{\mbox{\tiny[0]}}}^2 + 1 = 0. \]
Solving the above for initial velocity of the zeroth-order deformation equation does not make sense physically. In order to overcome this, we go back and modify homotopy $\mathcal{H}$ by the addition of an auxiliary operator \cite{BPL} $\Pi(\tilde{x};p)=\epsilon p(1-p)\cos(\tau)$ so that
\begin{multline*}
 \mathcal{H} \equiv (1-p)\bigg( \tilde{x}_{\tau\tau} + \lambda(p)^2 \tilde{x}\bigg) \\
-h(p) \bigg( \tilde{x}_{\tau\tau} + \lambda(p)^2\Big(\delta +\epsilon \cos\big(\lambda(1)\tau\big)\Big) \tilde{x} + \lambda(p)c\tilde{x}_{\tau}\bigg) + \epsilon p(1-p)\cos(\tau) =0.
\end{multline*}
This auxiliary operator contributes to resonant forcing, and hence corrects $\zeta^{\mbox{\tiny[{\em n}]}}$ at every deformation equation. It modifies the homotopy in a consistent manner, i.e.,  the auxiliary operator becoming zero at $p=0$ and at $p=1.$ The above homotopy works even if we wish to search for $4\pi$ periodic solutions. Continuing further, $x^{\mbox{\tiny [0]}}(\tau)$ remains same (Eqn.(\ref{zeroth_sol_damp})). Now, the first-order deformation equation is
\begin{multline}\label{first_order_damp}
 \ddot{x}^{\mbox{\tiny [1]}} + x^{\mbox{\tiny [1]}}=
\Big(- 2\lambda^{\mbox{\tiny [1]}}+ \left(-1+\delta + c\zeta^{\mbox{\tiny[0]}}\right) h^{\mbox{\tiny[1]}}
-\epsilon \Big) \cos(\tau) +\frac{\epsilon h^{\mbox{\tiny[1]}}}{2}\cos(2\tau)\\ - \Big( 2\zeta^{\mbox{\tiny[0]}} \lambda^{\mbox{\tiny[1]}}+ \left( c+
\zeta^{\mbox{\tiny[0]}}-\delta  \zeta^{\mbox{\tiny[0]}}\right) h^{\mbox{\tiny[1]}} \Big) \sin(\tau)+\frac{\epsilon \zeta^{\mbox{\tiny[0]}}h^{\mbox{\tiny[1]}}}{2}\sin(2\tau)+\frac{\epsilon h^{\mbox{\tiny[1]}}}{2}.
\end{multline}
The removal of the secular terms results in
\begin{gather*}
\zeta^{\mbox{\tiny[0]}} = \frac{1}{2c}\left(\frac {\epsilon \pm \sqrt {{{\epsilon}^{2} -4 {c}^{2}h^{\mbox{\tiny[1]}}}^{2}}
}{h^{\mbox{\tiny[1]}}}\right) \quad \mbox{and}
\quad	\lambda^{\mbox{\tiny[1]}} = \frac{1}{2}\bigg( - \epsilon  + \Big(-1 +\delta + c\zeta^{\mbox{\tiny[0]}}\Big)h^{\mbox{\tiny[1]}}\bigg).
\end{gather*}
We note that $\zeta^{\mbox{\tiny[0]}}$ is again a root of the  quadratic polynomial. Solving Eqn.(\ref{first_order_damp}) with the initial condition $(0,\zeta^{\mbox{\tiny [1]}}),$ we obtain
\begin{equation*}
x^{\mbox{\tiny [1]}}(\tau) = -\frac{\epsilon h^{\mbox{\tiny[1]}}}{3}\cos(\tau)-\frac{\epsilon h^{\mbox{\tiny[1]}}}{6}\cos(2\tau)
 + \left( \frac{\epsilon \zeta^{\mbox{\tiny[0]}} h^{\mbox{\tiny[1]}}}{3}+\zeta^{\mbox{\tiny[1]}} \right) \sin(\tau)
  -\frac{\epsilon\zeta^{\mbox{\tiny[0]}}h^{\mbox{\tiny[1]}}}{6} \sin(2\tau) 
  +\frac{\epsilon h^{\mbox{\tiny[1]}}}{2}.
\end{equation*}
We may proceed by substituting any of the two expressions for $\zeta^{\mbox{\tiny[0]}}$ in the above. Substituting the expression for $x^{\mbox{\tiny [0]}}(\tau)$ and 
$x^{\mbox{\tiny [1]}}(\tau)$ in the second-order deformation equation and eliminating the secular terms, we get 
\begin{multline*}
\zeta^{\mbox{\tiny[1]}} = \frac{1}{6(\epsilon - 2ch^{\mbox{\tiny[1]}}\zeta^{\mbox{\tiny[0]}})} \bigg( \Big(3{c}^{2}{\zeta^{\mbox{\tiny[0]}}}^{3}+ ( \delta+
2\epsilon-1) 3c{\zeta^{\mbox{\tiny[0]}}}^{2}+ ({c}^{2}+{\epsilon}^{2}) 3\zeta^{\mbox{\tiny[0]}}+3c\delta-2c\epsilon-3c \Big) {h^{\mbox{\tiny[1]}}}^{2} \\- \left( 3c{\zeta^{\mbox{\tiny[0]}}}^{2}+4\zeta^{\mbox{\tiny[0]}}{\epsilon}+3 c\right) \epsilon h^{\mbox{\tiny[1]}}+ ( {\zeta^{\mbox{\tiny[0]}}}^{2}+1 ) 3ch^{\mbox{\tiny[2]}}+6\epsilon\zeta^{\mbox{\tiny[0]}}\bigg) \quad \mbox{and}
\\
\lambda^{\mbox{\tiny[2]}} = -\frac{{\epsilon}^{2}}{4}+ \left(-1+ \delta- \frac{\epsilon}{2}-\frac{{\epsilon}^{2}}{3}-\frac{\delta \epsilon}{2} + c\zeta^{\mbox{\tiny[0]}}+c\zeta^{\mbox{\tiny[1]}} \right)h^{\mbox{\tiny[1]}}+ \frac{1}{2}\bigg(-1+ \delta+c\zeta^{\mbox{\tiny[0]}}\bigg) h^{\mbox{\tiny[2]}}\\ + \left(-\frac{1}{4} -\frac{\delta}{2} + {\frac {5{\epsilon}^{2}}{12}}+\frac{3{
\delta}^{2}}{4}+\frac{{c}^{2}{\zeta^{\mbox{\tiny[0]}}}^{2}}{4}+c\delta\zeta^{\mbox{\tiny[0]}}+\frac{2c\epsilon\zeta^{\mbox{\tiny[0]}}}{3}\right) {h^{\mbox{\tiny[1]}}}^{2}.
\end{multline*}

We continue with the same procedure upto $N$-th order deformation equation obtaining $\zeta^{\mbox{\tiny [{\em n}-1]}}$ and $\lambda^{\mbox{\tiny[{\em n}]}}$ by the removal of the secular terms at $n$-th order. At orders higher than one, $\zeta^{\mbox{\tiny[{\em n}]}}$ is obtained by solving a linear equation, resulting in a unique expression. As in the previous section, here too, we obtain $x^{\mbox{\tiny [{\em n}]}}(\tau),$ $\zeta^{\mbox{\tiny [{\em n}]}},$ $\lambda^{\mbox{\tiny[{\em n}]}}$ as functions of free parameters $h^{\mbox{\tiny[1]}},\cdots,h^{\mbox{\tiny[{\em n}]}},\;\delta$ and $\epsilon.$ We get an approximate expression for $\lambda(1)$ by substituting $\lambda^{\mbox{\tiny \tiny[1]}}$ to $\lambda^{\mbox{\tiny[{\em N}+1]}}$ in the truncated version of Eqn.(\ref{solution_taylor_lambda}) to obtain the relationship
\begin{equation}
	\lambda(1) = 1 = 1 + \sum_{n=1}^{N+1}\frac{1}{n!} \lambda^{\mbox{\tiny [{\em n}]}}.\label{taylor_approx1_damp}
\end{equation}
$2\pi$-periodic, the approximate solution to Eqn.(\ref{ode_lin_damp_math}), $x_N(t)$ can now be obtained, firstly by rescaling time $\tau=\tfrac{t}{\lambda(1)}=t$ in each $x^{\mbox{\tiny [{\em n}]}}(\tau),\;n=1,\cdots,N$ and secondly by using Eqn.(\ref{taylor_approx4}).

\begin{figure}[!ht]
\makebox[\textwidth][c]
{\includegraphics[draft=\status ,width=16cm,height=9.5cm]{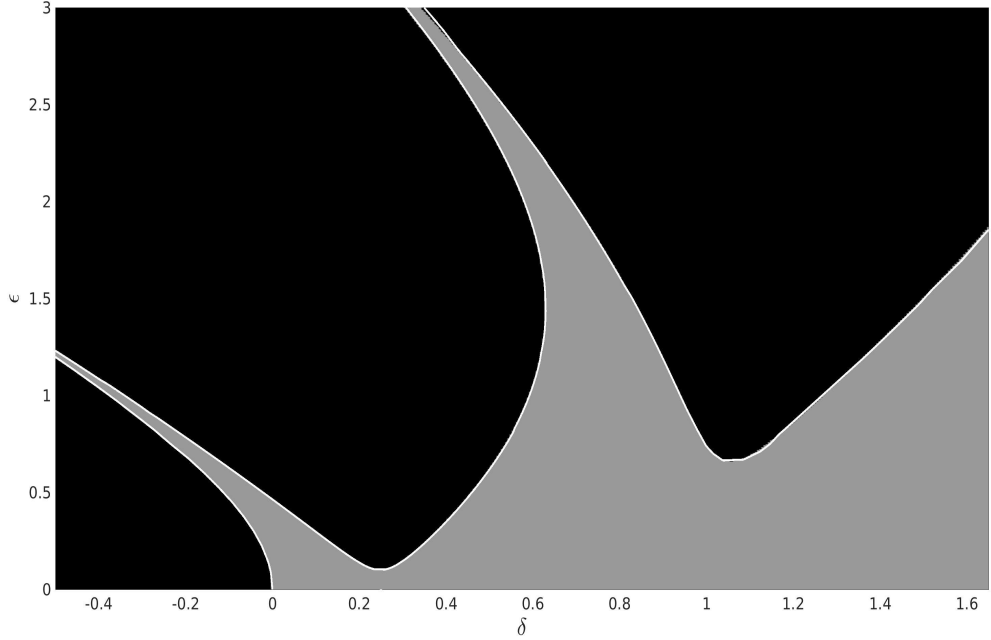}}
\caption{Stable regions (grey), unstable regions (black) and transition curves (white lines, homotopy$+$Galerkin) in the parameter plane with $c=0.1.$}
\label{Linear_Mathieu_Damp_Stab}
\end{figure}

\vspace{0.1in}
Equation(\ref{taylor_approx1_damp}) involves $N$+$3$ unknowns,  free parameters $h^{\mbox{\tiny[1]}}$ to $h^{\mbox{\tiny[{\em N}+1]}},$ $\delta$ and $\epsilon.$ As per Petrov-Galerkin scheme, we take the weighting functions same as basis functions of ${x}_N(t),$  i.e., $1,\ \cos(t),\  $ $\sin(t),\ \cos(2t),\ \sin(2t),\cdots$ and create $N+1$ equations by integrating the weighted residual over one time period of the solution,
$$
	\int_{0}^{2\pi} \cos(mt)\mathcal{R}_N(t) dt = 0,\quad
	\int_{0}^{2\pi} \sin(nt)\mathcal{R}_N(t) dt = 0,\quad
	\begin{rcases}
		\vspace{2mm}
			\begin{rcases} \vspace{1mm}
					m=0,1,\cdots,\tfrac{N+1}{2},\quad\cr
					n=1,2,\cdots,\tfrac{N-1}{2},\quad\cr		
			\end{rcases}
			\mbox{ for odd } N. 
		\cr 
			\begin{rcases} \vspace{1mm}
					m=0,1,\cdots,\tfrac{N}{2},\quad\cr
					n=1,2,\cdots,\tfrac{N}{2},\quad			
			\end{rcases}
			\mbox{ for even } N. 
	\end{rcases}
$$
We now have $N+2$ algebraic equations to be solved simultaneously for $N+3$ unknowns. Numerically solving this system of equations gives  transition curves in the parameter plane. Figure (\ref{Linear_Mathieu_Damp_Stab}) shows a region of $\delta$-$\epsilon$ plane covering $\delta \in [-0.5,1.7]$ and $\epsilon \in [0,3].$ Dividing it into $660\times900$ grid points and determining the stability with Floquet multipliers (Appendix \ref{floquet}) at each grid point, the stable and unstable regions are generated. Also generated are the transition curves (white lines in Fig.(\ref{Linear_Mathieu_Damp_Stab})) obtained numerically by solving a system of six equations ($N=4$) using Maple. While solving the algebraic equations, we keep the precision of the numerical solver high ($30$ digits) to ensure accurate computation of $\delta.$ We observe that the root of $\zeta^{\mbox{\tiny[0]}}$ with negative discriminant enables determination of the entire transition curve corresponding to a given unstable region while the other root fails to do so.

\section{Mathieu equation with impulsive excitation}
\label{sec:4}
We now consider the linear Mathieu equation with impulsive parametric excitation periodic with period $2\pi$
\begin{equation} \label{ode_lin_impulse_math}
	 \ddot{x} + \bigg(\Delta  + \epsilon \sum_{k=0}^{\infty}\delta\Big(t-(2k+1)\pi \Big)\bigg) x = 0. 
\end{equation}
By $\delta,$ we denote Dirac delta function, and introduce $\Delta$ and $\epsilon$ as excitation parameters. The stability analysis of solutions of oscillators with periodic impulsive parametric excitation is studied in depth \cite{impulse}. We retain the linear operator to be a simple harmonic oscillator and choose nonlinear operator to be Mathieu equation with impulsive parametric excitation. Let homotopy be
\begin{equation} \label{PDE_lin_impulse_unscaled}
	\mathcal{H} \equiv (1-p)\left(\tilde{x}_{tt} + \tilde{x}\right) - h(p) \Bigg( \tilde{x}_{tt} + \bigg(\Delta  + \epsilon \sum_{k=0}^{\infty}\delta\Big(t-(2k+1)\pi \Big)\bigg) \tilde{x}\Bigg)=0.
\end{equation}
\vspace{0.1in}

By scaling time $t$ to $\tau=\tfrac{t}{\lambda(p)}$ in Eqn.(\ref{PDE_lin_impulse_unscaled}) and replacing the forcing part $\delta\Big(\lambda(p)\tau-(2k+1)\pi \Big)$ by $\delta\Big(\lambda(1)\tau-(2k+1)\pi \Big),$ we arrive at
\begin{equation*}
	\mathcal{H} \equiv (1-p)\bigg( \tilde{x}_{\tau\tau} + \lambda(p)^2 \tilde{x}\bigg)  
-h(p) \Bigg( \tilde{x}_{\tau\tau} + \lambda(p)^2\bigg(\Delta  + \epsilon \sum_{k=0}^{\infty}\delta\Big(\lambda(1)\tau-(2k+1)\pi \Big)\bigg)\tilde{x} \Bigg)=0. 
\end{equation*}
We begin by searching for $2\pi$ periodic solution ($\lambda(1)=1$) with initial condition (1,0). As usual, $\lambda(0)=1.$ Compared to the previous treatment of homotopy analysis, for impulsive parametric excitation, we differ in an important aspect in the process of the removal of secular terms, mainly for the sake of computational convenience. As we are looking for a periodic solution, it suffices to construct only one period of the solution computationally. Therefore, we restrict the homotopy over, say, first period of the solution.
This leads to a modification of the homotopy as 
\begin{equation}\label{PDE_lin_impulse_scaled_lambda1_final}
	\mathcal{H}_{[0,2\pi]} \equiv (1-p)\bigg( \tilde{x}_{\tau\tau} + \lambda(p)^2 \tilde{x}\bigg)  
-h(p) \bigg( \tilde{x}_{\tau\tau} + \lambda(p)^2\Big(\Delta  + \epsilon \delta\big(\tau-\pi \big)\Big)\tilde{x} \bigg)=0. 
\end{equation}
Subsequent analysis in this section is restricted to $\tau \in [0,2\pi].$ Let $\hat{x}_N(t)$ be the approximate solution to the proposed homotopy at $p=1.$ The zeroth-order deformation equation for Eqn.(\ref{PDE_lin_impulse_scaled_lambda1_final}) is
\begin{equation*}
 \mathcal{L}(x^{\mbox{\tiny [0]}}) = \ddot{x}^{\mbox{\tiny [0]}} + x^{\mbox{\tiny [0]}}=0.
\end{equation*}
We choose the initial conditions for various order deformation equations here in the same manner as for the linear Mathieu equation (Eqn.(\ref{ode_lin_math})), i.e.,
\begin{equation}\label{initial_conditions_nth_order_impulse}
x^{\mbox{\tiny [0]}}(0)=1, \quad \dot{x}^{\mbox{\tiny [0]}}(0)=0\quad \mbox{and}	\quad x^{\mbox{\tiny[{\em n}]}}(0)=0, \quad \dot{x}^{\mbox{\tiny[{\em n}]}}(0)=0, \quad n\geq1. 
\end{equation}
Accordingly, we obtain at zeroth-order, 
\begin{equation*}
	x^{\mbox{\tiny [0]}}(\tau) = \cos(\tau).
\end{equation*}
Substituting the above in first-order deformation equation, we get
\begin{equation}\label{blahblah}
\ddot{x}^{\mbox{\tiny [1]}} + x^{\mbox{\tiny [1]}}
 = \bigg( -2 \lambda^{\mbox{\tiny [1]}} -\Big( 1-\Delta - \epsilon\delta (\tau-\pi)\Big)h^{\mbox{\tiny [1]}} \bigg) \cos(\tau)=F_{1}(\tau,\lambda^{\mbox{\tiny[1]}}).
\end{equation}
Coefficient of $\cos(\tau)$ is not a constant, but depends on $\tau.$ So by removing the resonant forcing using
\begin{gather*}
\int_{0}^{2\pi}\cos(\tau)F_{1}(\tau,\lambda^{\mbox{\tiny[1]}}) d\tau = 0,
\end{gather*}
gives
\begin{equation}\label{lambda_impulse}
\lambda^{\mbox{\tiny [1]}} =\frac{h^{\mbox{\tiny [1]}}}{2} \left(-1+ \Delta+\frac{\epsilon}{\pi} \right). 
\end{equation}
Inserting Eqn.(\ref{lambda_impulse}) in Eqn.(\ref{blahblah}), we get
\begin{equation*}\label{first_order_impulse_1}
\ddot{x}^{\mbox{\tiny [1]}} + x^{\mbox{\tiny [1]}}
 = \epsilon h^{\mbox{\tiny [1]}} \bigg(-\frac{1}{\pi} + \delta (\tau-\pi)  \bigg)\cos(\tau).
\end{equation*}
By solving the above with the appropriate initial condition, we get
\begin{equation*}
x^{\mbox{\tiny [1]}}(\tau) = \frac{\epsilon h^{\mbox{\tiny [1]}}}{2}\left( 2 H\left(\tau-\pi \right) - {\frac {\tau}{\pi}} \right)\sin(\tau),
\end{equation*}
where $H$ denotes Heaviside function. Repeating the above steps upto $N$-th order deformation equation, we obtain $x^{\mbox{\tiny [{\em n}]}}(\tau),$ $\zeta^{\mbox{\tiny [{\em n}]}}$ and $\lambda^{\mbox{\tiny[{\em n}]}}$ as functions of parameters $h^{\mbox{\tiny[1]}},\cdots,h^{\mbox{\tiny[{\em n}]}},\;\Delta$ and $\epsilon.$ Substituting the expressions for $\lambda^{\mbox{\tiny \tiny[1]}}$ to $\lambda^{\mbox{\tiny[{\em N}+1]}}$ in the equation below, we get a relationship between $\delta$ and $\epsilon$
\begin{equation}\label{taylor_approx1_impulse}
	\lambda(1) = 1 = 1 + \sum_{n=1}^{N+1}\frac{1}{n!} \lambda^{\mbox{\tiny [{\em n}]}}.
\end{equation}
For $N=3,$ this expression is given in Appendix \ref{lin_impulse_eqns}. $\hat{x}_N(t),$ an approximate solution to Eqn.(\ref{ode_lin_impulse_math}) is obtained by scaling time $\tau=\tfrac{t}{\lambda(1)}=t$ in each $x^{\mbox{\tiny [{\em n}]}}(\tau),\;n=1,\cdots,N$ and  using Eqn.(\ref{taylor_approx4}). For $N=3,$ we have 
\begin{dmath}\label{x_sol_1}
\hat{x}_3(t)=
\cos(t) 
-\frac{\epsilon^{2} h^{\mbox{\tiny[1]}}}{2} \bigg( 3{h^{\mbox{\tiny[1]}}} + h^{\mbox{\tiny[2]}} + {\frac {\left( 2 \pi \Delta+\epsilon \right) {h^{\mbox{\tiny[1]}}}^{2}}{\pi}} \bigg)\cos(t) H(t-\pi)
+\epsilon \bigg(3 h^{\mbox{\tiny[1]}}+ h^{\mbox{\tiny[2]}}+\frac{\,h^{\mbox{\tiny[3]}}}{6} + {\frac {  \left(6\,{\pi} \Delta+ 3 \epsilon \right) {h^{\mbox{\tiny[1]}}}^{2}}{2{\pi} }}+{\frac { \left( 2{\pi} \Delta+ \epsilon \right) h^{\mbox{\tiny[1]}}h^{\mbox{\tiny[2]}}}{12{\pi} }} +{\frac { \Big( (1- {\pi}^{2}){\epsilon}^{2}+4\pi\,\Delta( \pi \Delta + \epsilon) \Big) {h^{\mbox{\tiny[1]}}}^{3}}{4{\pi}^{2}}}\bigg)\sin(t) H(t-\pi)
-\frac{\epsilon}{2\pi}\bigg({{3h^{\mbox{\tiny[1]}}}}+{{h^{\mbox{\tiny[2]}}}}+{ \frac{h^{\mbox{\tiny[3]}}}{6}}+{\frac{ 3\left( 2{\pi}\Delta +  \epsilon \right) {h^{\mbox{\tiny[1]}}}^{2}}{2{\pi}}}+{\frac {\left( 2{\pi} \Delta+ \epsilon \right) h^{\mbox{\tiny[1]}}h^{\mbox{\tiny[2]}}}{2{\pi}}}+{\frac {\Big( 12{\pi}{\Delta}(\pi \Delta+ \epsilon) +(3-\pi^2){\epsilon}^{2} \Big) {h^{\mbox{\tiny[1]}}}^{3}}{12{\pi}^{2}}}\bigg) t\sin(t)
+\epsilon^2 \bigg({\frac { \left( 2{\pi}\Delta +{\epsilon} \right) {h^{\mbox{\tiny[1]}}}^{3}}{2{\pi}^{2}}}+{\frac {3{h^{\mbox{\tiny[1]}}}^{2}}{2\pi}}+{\frac {h^{\mbox{\tiny[1]}}h^{\mbox{\tiny[2]}}}{2\pi}}\bigg) t\cos(t) H(t-\pi)
+\frac {{\epsilon}^{3}{h^{\mbox{\tiny[1]}}}^{3}}{4{\pi}} t \sin(t) H(t-\pi)
-\frac{{\epsilon}^{2}h^{\mbox{\tiny[1]}}}{8\pi^2}\bigg({3{h^{\mbox{\tiny[1]}}}}+{h^{\mbox{\tiny[2]}}} + {\frac { \left( 2\pi\,\Delta + {\epsilon} \right) {h^{\mbox{\tiny[1]}}}^{2}}{{\pi}}} \bigg)t^2\cos(t)
- {\frac {{\epsilon}^{3}{h^{\mbox{\tiny[1]}}}^{3}}{8{\pi}^{2}}} t^2\sin(t) H(t-\pi)
+ {\frac {{\epsilon}^{3}{h^{\mbox{\tiny[1]}}}^{3}}{48{\pi}^{3}}} t^3\sin(t) .
\end{dmath}
Although $x^{\mbox{\tiny [{\em n}]}}(\tau)$ has polynomially growing terms, the required validity of Eqn(\ref{x_sol_1}) is over $t\in [0,2\pi]$ and over this bounded domain, $\hat{x}_3(t)$ is bounded.

\begin{figure}
  \begin{subfigure}[b]{0.49\textwidth}
  \centering
    \includegraphics[draft=\status ,height=4cm,width=7cm]{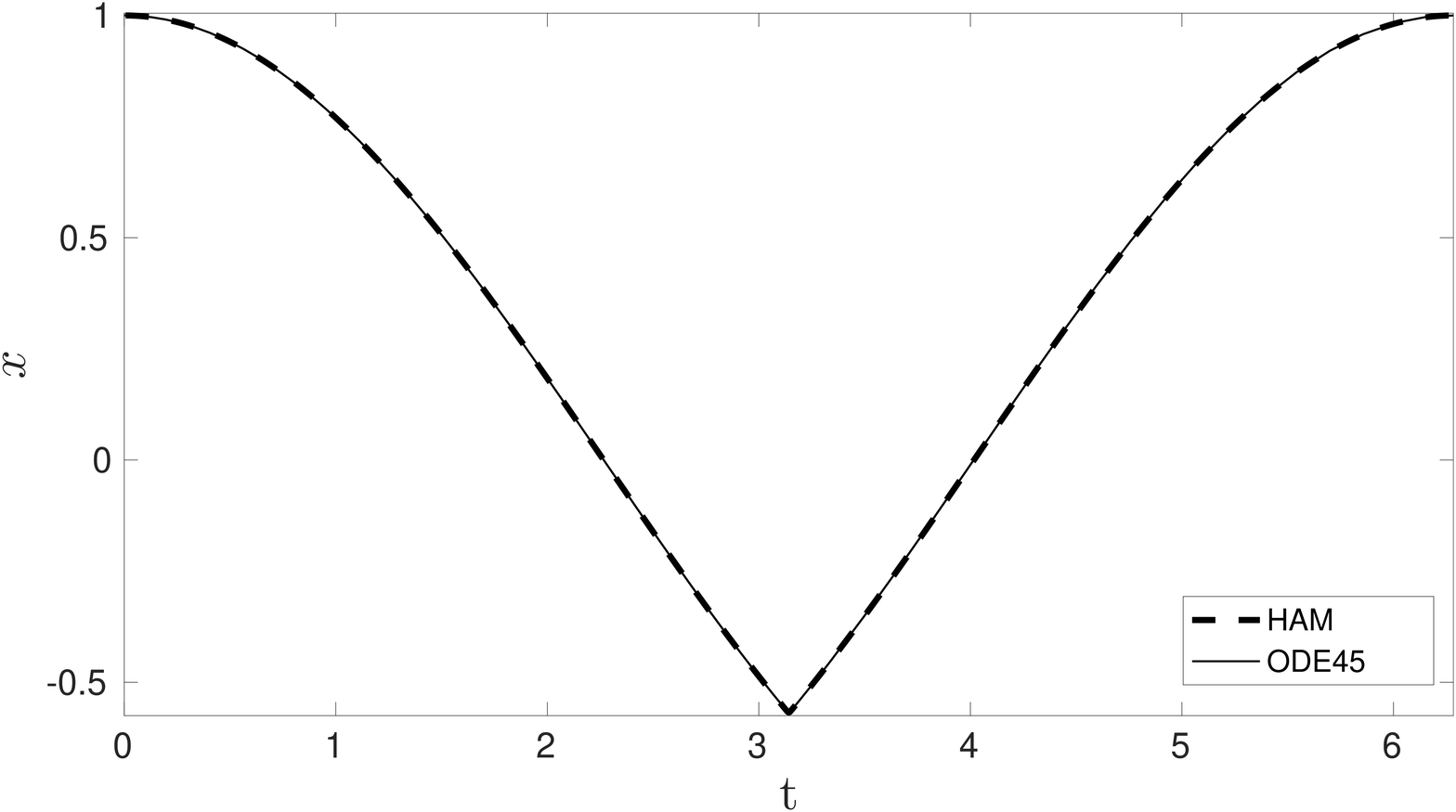}
  \end{subfigure}
%  \hfill
  \begin{subfigure}[b]{0.49\textwidth}
\centering
    \includegraphics[draft=\status ,height=4cm,width=7cm]{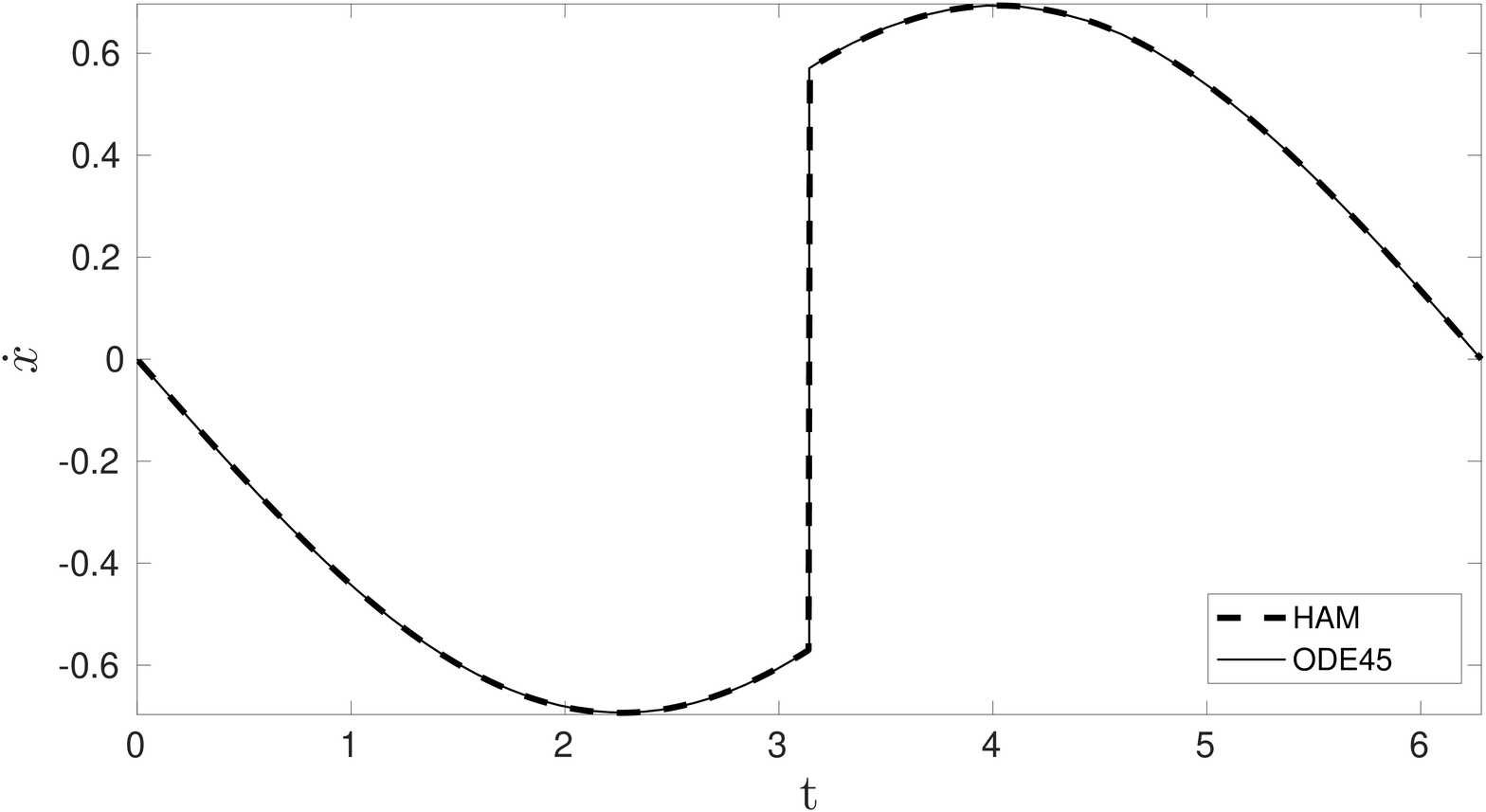}
  \end{subfigure}
  \caption{Comparison of $2\pi$ periodic solution for $\Delta = 0.4802..,\  \epsilon=2.0,\ N=3.$}
  \label{solution_compare_1}
\end{figure}

Equation (\ref{taylor_approx1_impulse}) involves unknowns $h^{\mbox{\tiny[1]}},\cdots,h^{\mbox{\tiny[{\em N}+1]}}.$ As earlier, we apply Petrov-Galerkin scheme, taking weighting functions $t,\ \cos(t), \ \cos(t)H(t-\pi), \ \sin(t)H(t-\pi),\cdots$ to create $N+1$ equations by integrating the weighted residual over time period $2\pi.$ Then we have $N+2$ algebraic equations to be solved for $N+3$ unknowns. For a particular case of $N=3$ and $\epsilon=2,$ we solve the system of equations (Appendix \ref{lin_impulse_eqns}) numerically for free parameters and $\Delta.$ By substituting the values so obtained in Eqn.(\ref{x_sol_1}), we plot the solution over the first period and see from Fig(\ref{solution_compare_1}) that it mimics the solution obtained by numerical integration of Eqn.(\ref{ode_lin_impulse_math}).

\vspace{0.1in}

Now we choose $\lambda(1)=2$ to look for $4\pi$ periodic solution with initial condition $(0,1)$. Homotopy $\mathcal{H}$ is modified to
\begin{equation*}
\mathcal{H}_{[0,4\pi]} \equiv (1-p)\bigg( \tilde{x}_{\tau\tau} + \lambda(p)^2 \tilde{x}\bigg)  
-h(p) \Bigg( \tilde{x}_{\tau\tau} + \lambda(p)^2\bigg(\Delta  + \frac{\epsilon }{2}\left(\delta\Big(\tau-\frac{\pi}{2}\Big) + \delta\Big(\tau-\frac{3\pi}{2}\Big)\right)\bigg)\tilde{x} \Bigg)=0.
\end{equation*}
It has two terms from impulsive parametric excitation that act over interval zero to $4\pi.$ The zeroth-order deformation equation for the above homotopy would be same as before, i.e.,
\begin{equation*}
 \mathcal{L}(x^{\mbox{\tiny [0]}}) = \ddot{x}^{\mbox{\tiny [0]}} + x^{\mbox{\tiny [0]}}=0.
\end{equation*}
The initial condition for the zeroth-order deformation equation is $x^{\mbox{\tiny [0]}}(0)=0$ and $\dot{x}^{\mbox{\tiny [0]}}(0)=1$ resulting in 
\begin{equation}\label{zeroth_order_solution_impulse_4} 
	x^{\mbox{\tiny [0]}}(\tau) = \sin(\tau).
\end{equation}
The initial conditions for higher order deformation equations are fixed in accordance with Eqn.(\ref{initial_conditions_nth_order_impulse}). By substituting Eqn.(\ref{zeroth_order_solution_impulse_4}) in  first-order deformation equation, we get
\begin{equation*}
\ddot{x}^{\mbox{\tiny [1]}} + x^{\mbox{\tiny [1]}}
 = \bigg(-2\lambda^{\mbox{\tiny [1]}}-\left(1-\Delta-\frac{\epsilon}{2}\delta\Big(\tau-\frac{\pi}{2} \Big)-\frac{\epsilon}{2}\delta\Big(\tau-\frac{3\pi}{2} \Big)\right) h^{\mbox{\tiny [1]}}\bigg) \sin(\tau)=F_{1}(\tau,\lambda^{\mbox{\tiny[1]}}).
\end{equation*}
By eliminating the secular terms using
\begin{gather*}
\int_{0}^{4\pi}\sin(\tau)F_{1}(\tau,\lambda^{\mbox{\tiny[1]}}) d\tau = 0
\end{gather*}
and then solving for $\lambda^{\mbox{\tiny [1]}},$ we get
\begin{equation}\label{lambda_impulse_4}
\lambda^{\mbox{\tiny [1]}} =\frac{h^{\mbox{\tiny [1]}}}{2} \left(-1+ \Delta+\frac{\epsilon}{\pi} \right) .
\end{equation}
Substituting Eqn.(\ref{lambda_impulse_4}) in the first-order equation, we get
\begin{equation*}
\ddot{x}^{\mbox{\tiny [1]}} + x^{\mbox{\tiny [1]}}
 = \frac{\epsilon h^{\mbox{\tiny [1]}}}{2}\left(-\frac{2}{\pi}+\delta\Big(\tau-\frac{\pi}{2} \Big)+\delta\Big(\tau-\frac{3\pi}{2} \Big)\right) \sin(\tau).
\end{equation*}
By solving the above equation, we get
\begin{equation*}
x^{\mbox{\tiny [1]}}(\tau) = \frac{\epsilon h^{\mbox{\tiny [1]}}}{2} \bigg( - H \left(\tau-\frac{3\pi}{2} \right) - H\left(\tau-\frac{\pi}{2} \right)+\frac{\tau}{\pi}\bigg) \cos(\tau)-\frac{\epsilon{h^{\mbox{\tiny [1]}}}}{2\pi}\sin(\tau).
\end{equation*}
We repeat the same upto the $N$-th order deformation equation. Again, by approximating $\lambda(1)$ using expressions for $\lambda^{\mbox{\tiny \tiny[1]}}$ to $\lambda^{\mbox{\tiny[{\em N}+1]}}$ as obtained via removal of the secular terms, we get the relationship
\begin{equation}\label{lambda_taylor_impulse_2}
	\lambda(1) = 2 = 1 + \sum_{n=1}^{N+1}\frac{1}{n!} \lambda^{\mbox{\tiny [{\em n}]}}.
\end{equation}
$\hat{x}_N(t),$ the solution to Eqn.(\ref{ode_lin_impulse_math}) over $t \in [0,4\pi]$ is approximated by rescaling time $\tau=\tfrac{t}{\lambda(1)}=\tfrac{t}{2}$ in Eqn.(\ref{taylor_approx4}). We compute up to order $N=3$ for obtaining $\hat{x}_3(t)$ (Eqn.(\ref{x_sol_2}), Appendix \ref{lin_impulse_sol}).

\begin{figure}[!h]
    \begin{subfigure}[b]{0.49\textwidth}
  \centering
    \includegraphics[draft=\status ,height=4cm,width=7cm]{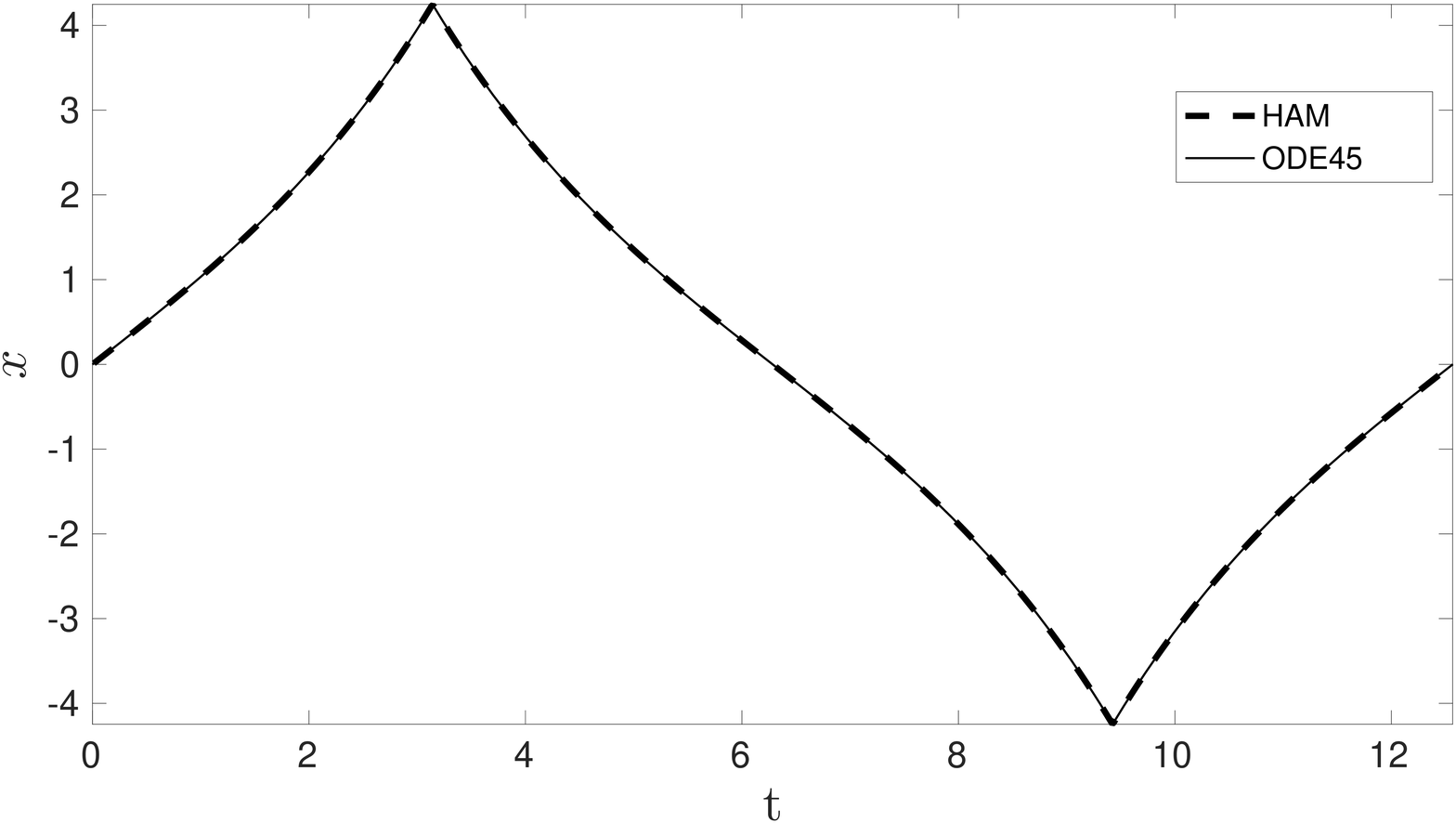}
  \end{subfigure}
%  \hfill
  \begin{subfigure}[b]{0.49\textwidth}
\centering
    \includegraphics[draft=\status,height=4cm,width=7cm]{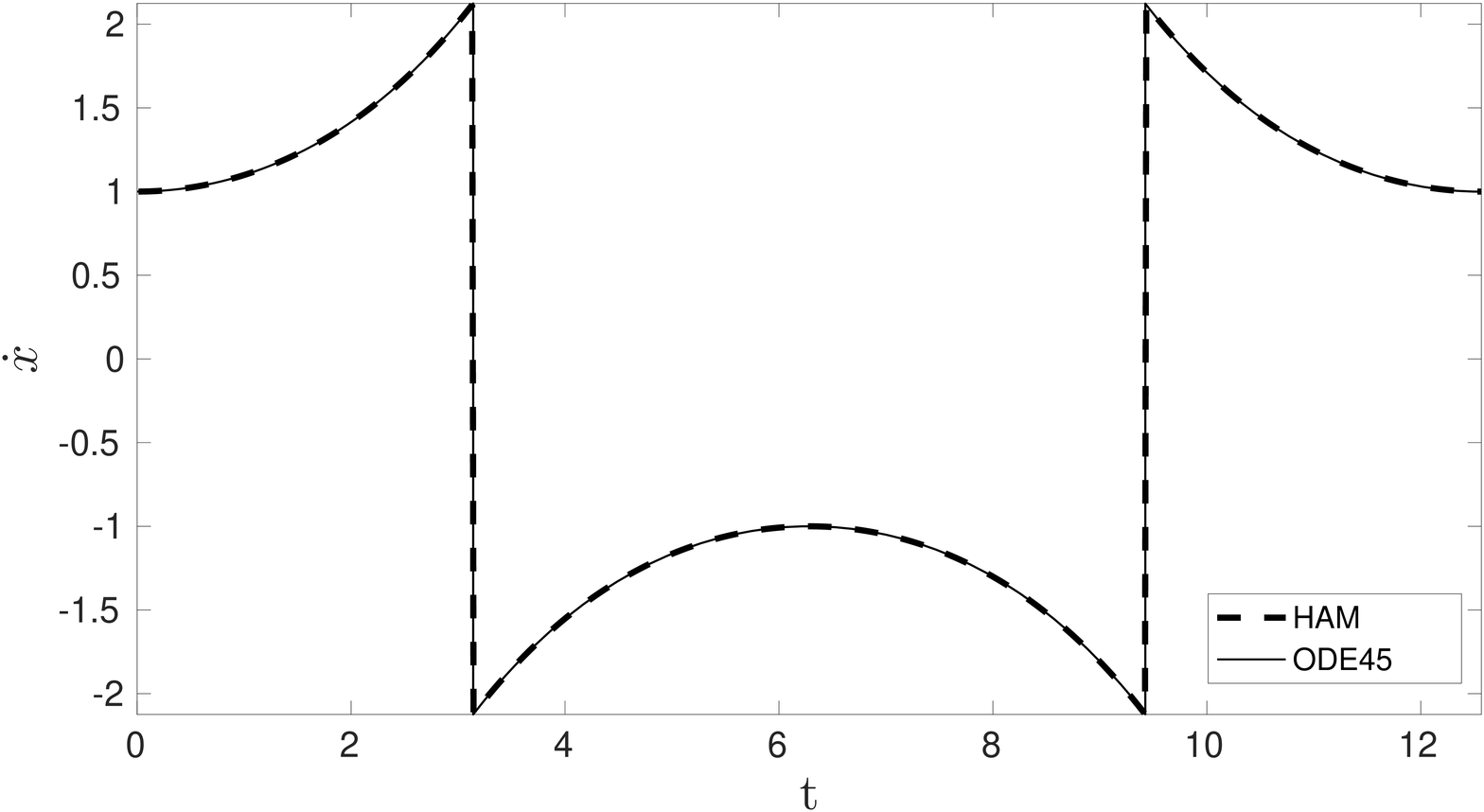}
  \end{subfigure}
      \caption{Comparison of $4\pi$ periodic solution for $\Delta = -0.1945..,\  \epsilon=1.0,\ N=3.$}
      \label{solution_compare_2}
\end{figure}
\begin{figure}[!h]
\makebox[\textwidth][c]
{\includegraphics[draft=\status ,width=18cm,height=10cm]{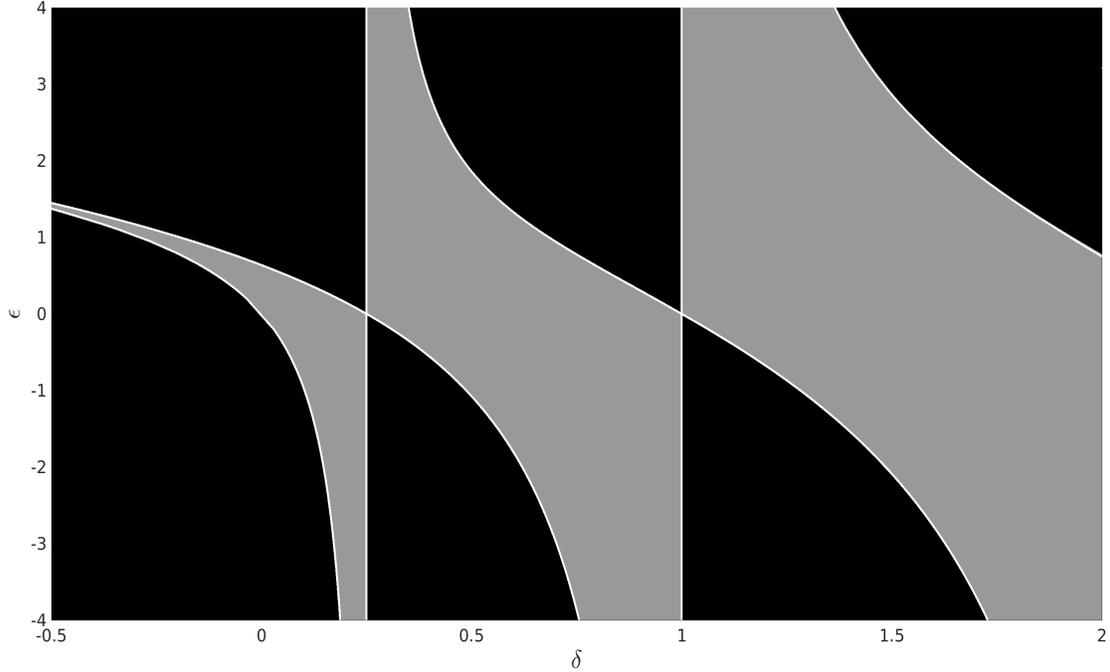}}
\caption{Stable regions (grey), unstable regions (black) and transition curves (white lines, homotopy$+$Galerkin) in the parameter plane.}
\label{Linear_Mathieu_impulse_Stab}
\end{figure}

Following the Petrov-Galerkin scheme, we take $t,\ t\cos(t), \ H(t-\pi)\cos(t), \ H(t-\pi)\sin(t),\ H(t-3\pi)\cos(t), \ H(t-3\pi)\sin(t),\cdots$ as weighting functions and  create $N+1$ equations by integrating the weighted residual over one time period, $4\pi.$ For a fixed value of $\epsilon,$ say $1$, Eqn.(\ref{lambda_taylor_impulse_2}) along with $N+1$ equations via Galerkin projections may be solved numerically for $N+1$ unknowns and $\Delta$. By substituting values so obtained in Eqn.(\ref{x_sol_2})(Appendix \ref{lin_impulse_sol}) and plotting the solution over the first period in Fig.(\ref{solution_compare_2}), we find that it matches well with the solution obtained by ode-45 integration. The stability regions in the parameter plane are shown in Fig.(\ref{Linear_Mathieu_impulse_Stab}). The region, with $\Delta$ varying from $-0.5$ to $2$ and $\epsilon$ varying from $-4$ to $4$ is divided into $780\times2800$ grid points. As previously done, the stability regions are generated by calculating the Floquet multipliers (Appendix \ref{floquet}) for each grid point. As can be seen from Fig.(\ref{Linear_Mathieu_impulse_Stab}), the right branch of each unstable region in $\epsilon>0$ half-plane and the left branch for the same in $\epsilon<0$ half-plane is a straight line. A periodic solution along any branch that is a straight line, $\hat{x}_N(t)$ is independent of $N$ and is either $\sin(t)$ or $\cos(\tfrac{t}{2})$ exactly. The algebraic equations we get in this case via Galerkin projections are either trivially zero or $\Delta=1$ or $\Delta=0.25$ depending on whether period is $2\pi$ or $4\pi$ respectively. Also, the $\Delta$-$\epsilon$ relationship via homotopy plays no role in determining the transition curves. The remaining transition curves are obtained by routinely solving the system of algebraic equations ($2\pi$-periodic case with $N=3$, Appendix \ref{lin_impulse_eqns}) and plotted in Fig.(\ref{Linear_Mathieu_impulse_Stab})(white lines).

%%%%%%%%%%%%%%%%%%%%%%%%%%%%%%%%%%%%%%%%%%%%%%%%%%%%%%
\section{Conclusion and further work}

We have combined the homotopy analysis method with the Galerkin projections and applied it in the context of parametrically excited oscillators. This primarily analytical attempt has enabled us to approximate not only periodic solutions on the stability oundaries but also  transition curves in the parameter plane. As parameter $\epsilon$ need not be assumed to be small in the  framework of homotopy analysis, with the help of this technique we may find transition curves spread over large regions of the parameter plane with affordable computation. We demonstrate the results obtained via this method with three related, yet distinct versions of linear Mathieu equation by exploiting the flexibility inherent to homotopy analysis method. The  method is as analytical as the method of harmonic balance.

The method can easily be employed to find transition curves of a variety of Mathieu-like equations, for example, scalable Mathieu-like oscillators, Mathieu equation with delay, fractional-order Mathieu and similar versions of Ince's equation etc. We also believe that the method will find application in approximating not only natural, periodic responses but forced responses of strongly nonlinear oscillators also. It may also find application in getting approximate solution of nonlinear boundary value problems.
\\

\section*{Acknowledgements}
We sincerely thank Prof. Anindya Chatterjee for providing technical comments on this article. We also thank Dr. Hari Nair for proofreading this article.
\\~\\
\textbf{Funding}: This research did not receive any specific grant from funding agencies in the public, commercial, or not-for-profit sectors.

%%%%%%%%%%%%%%%%%%%%%%%%%%%%%%%%%%%%%%%%%%%%%%%%%%%%%%
\iffalse
\appendix
\section*{Appendix}
\fi

\appendixpage
\begin{appendices}

\section{Numerical calculation of Floquet multipliers}
\label{floquet}
Consider linear Mathieu equation in the form
\begin{equation}    \label{eq:main}
    \dot{x}=A(t)x \quad \mbox{with} \quad
 A(t)=\left[ \begin{matrix}
   0 & 1  \\
   -\Big(\delta +\epsilon \cos (t)\Big) & 0  \\
\end{matrix} \right].
\end{equation}
$A(t)$ is periodic in time with period $2\pi.$ Let $\mathrm{x}_1(t)$  and $\mathrm{x}_2(t)$ be solutions to Eqn.(\ref{eq:main}) with the initial conditions $(1,0)$ and $(0,1)$ respectively. To obtain the two solutions at a given point in the parameter plane (corresponding $\delta$ and $\epsilon$ values) for a period $2\pi,$ we integrate Eqn.(\ref{eq:main}) numerically using ode45, built-in integrator of MATLAB. The Floquet matrix is then obtained by
\[\tilde{X}=\left[ \begin{matrix}
   {\mathrm{x}_{1}}(T) & {\mathrm{x}_{2}}(T)  
\end{matrix} \right].
\]
Eigenvalues $\lambda_{1}$ and $\lambda_{2}$ of $\tilde{X}$ are nothing but the Floquet multipliers. If $\left| {{\lambda }_{1}} \right|>1$ or $\left| {{\lambda }_{2}} \right|>1$, then the solution in general grows exponentially and the system is unstable. The system is stable if $\left| {{\lambda }_{i}} \right|=1$ for $i=1,2.$ The point is marked as black in the parameter plane when the system is unstable;  or else, the point is marked grey (stable). By repeating this over all the grid points, we generate the stability chart for the equation. For damped the Mathieu equation, we have
\[A(t)=\left[ \begin{matrix}
   0 & 1  \\
   -\Big(\delta +\epsilon \cos (t)\Big) & -c  \\
\end{matrix} \right].\]
We arrive at the Floquet matrix following the same procedure as for the undamped case. The system is unstable if the same condition holds as in the case of the undamped Mathieu. The point is marked as black in the parameter plane if the system is unstable; or else, the point is marked grey (stable). By repeating this over all the grid points, the stability chart for the equation is generated. Note that in this case, the system is stable provided $\left| {{\lambda }_{i}} \right|<1$ for both $i=1,2.$

\vspace{0.1in}

The impulsive Mathieu equation (periodic in $t\in[0,2\pi]$) given by
\begin{equation}  \label{eq:impulse_short}
	 \ddot{x} + \Big(\Delta  + \epsilon \delta\big(t-\pi \big)\Big) x = 0,
\end{equation}
when written in the form of Eqn.(\ref{eq:main}) gives
\[A(t)=\left[ \begin{matrix}
   0 &  & 1  \\
  -\Big(\Delta  + \epsilon \delta\big(t-\pi \big)\Big) & & 0  \\
\end{matrix} \right].\]
Here $\delta$ symbol denotes Dirac delta function. To get the Floquet multipliers, we integrate Eqn.(\ref{eq:impulse_short}) with initial conditions $(1,0)$ and $(0,1)$ over $t \in [0,2\pi].$ Since the first impulse acts at $t=\pi,$ causing a jump in the velocity, we  integrate over $t\in [0,2\pi]$ in two parts, first over $t\in [0,\pi]$ and second $t\in[\pi,2\pi].$ We choose the initial condition for the first part either $(1,0)$ or $(0,1)$ and for the second part according to the following. By denoting the time instant just after impulse by $t=\pi+\mu$ ($\mu$ being an infinitesimally small quantity), we integrate both sides of the Eqn.(\ref{eq:impulse_short}) from $t=0$ to $\pi+\mu$
\begin{equation*}
\int\limits_{0}^{\pi +\mu }{\ddot{x} \,dt}+\Delta \int\limits_{0}^{\pi +\mu }{x \,dt}+\epsilon \int\limits_{0}^{\pi +\mu }{\delta(t-\pi )}x\ dt =0 \quad \mbox{to get}
\end{equation*}
\begin{equation}  \label{eq:impulseVelo}
\dot{x}(\pi +\mu )=\dot{x}(0)-\epsilon x(\pi)-\Delta \int\limits_{0}^{\pi +\mu }{x \,dt}.
\end{equation}
The last term on the right hand side of Eqn.(\ref{eq:impulseVelo}) is to be obtained by computing the area under (the numerically sought) $x(t)$ using the trapezoidal method. Using $x(\pi+\mu)=x(\pi)$ and Eqn.(\ref{eq:impulseVelo}), we get the state vector just after the impulse. With this state vector as the initial condition, we integrate Eqn.(\ref{eq:impulse_short}) from $t=\pi+\mu$ to $2\pi$ numerically, again using ode45 integrator to get the state at $t=2\pi.$ Then we get the Floquet matrix and its eigenvalues. The procedure discussed for the Mathieu equation without damping is carried out for all the grid points in the parameter plane for obtaining the stability chart.

\section{Algebraic equations via homotopy and Galerkin projections}
\subsection{Linear Mathieu equation}
\label{lin_math_eqns}
	
Considering Eqn.(\ref{taylor_approx1}) with $N=3$ and substituting expressions for $\lambda^{\tiny [1]},\ \lambda^{\tiny [2]},\ \lambda^{\tiny [3]}$ and $\lambda^{\tiny [4]}$ gives
\begin{dmath*}
1+\left( 1- \delta-\frac{\epsilon}{2} \right) 2h^{\mbox{\tiny [1]}}+ \left( 1-\delta-\frac{\epsilon}{2} \right) \frac{3 h^{\mbox{\tiny [2]}}}{4}+ \left( 1-\delta-\frac{\epsilon}{2} \right) \frac{h^{\mbox{\tiny [3]}}}{6}+ \left( {1}-\delta-{\frac {\epsilon}{2}} \right) \frac{h^{\mbox{\tiny [4]}}}{48}+ \left(\frac{1}{2}+\delta+\frac{\epsilon}{2}-\frac{3}{2}-\frac{3 \epsilon \delta}{2}-{\frac {5 {\epsilon}^{2}}{16}} \right)\frac{3 {h^{\mbox{\tiny [1]}}}^{2}}{2}+ \left( \frac{1}{2}+\delta+\frac{\epsilon}{2}-{\frac {3 {\delta}^{2}}{2}}-{\frac {3 \epsilon \delta}{2}}-{\frac {5 {\epsilon}^{2}}{16}} \right) \frac{3 h^{\mbox{\tiny [1]}}h^{\mbox{\tiny [2]}}}{4}+ \left( \frac{1}{2}+{\delta}+\frac{\epsilon}{2}-{\frac {3 {\delta}^{2}}{2}}-{\frac {3 \epsilon \delta}{2}}-{\frac {5 {\epsilon}^{2}}{16}} \right) \frac{{h^{\mbox{\tiny [2]}}}^{2}}{16}+ \left( \frac{1}{6}+\frac{\delta}{3}+\frac{\epsilon}{6}-\frac{{\delta}^{2}}{2}-\frac{ \epsilon \delta}{2}-{\frac {5 {\epsilon}^{2}}{48}} \right) \frac{ h^{\mbox{\tiny [1]}}h^{\mbox{\tiny [3]}}}{4}+ \left( 1+{\delta}+\frac{\epsilon}{2}+3 {\delta}^{2}+3 \epsilon \delta+{\frac {5 {\epsilon}^{2}}{8}}-5 {\delta}^{3}-{\frac {15 \epsilon {\delta}^{2}}{2}}-{\frac {25 {\epsilon}^{2} \delta}{8}}-{\frac {19 {\epsilon}^{3}}{64}} \right) \frac{{h^{\mbox{\tiny [1]}}}^{3}}{4}+ \left( {1+{\delta}+{\frac {\epsilon}{2}}+ {3 {\delta}^{2}}+3 \epsilon \delta+{\frac {5 {\epsilon}^{2}}{8}}-5 {\delta}^{3}-\frac {15 \epsilon {\delta}^{2}}{2}-\frac {25 {\epsilon}^{2} \delta}{8}-\frac {19 {\epsilon}^{3}}{64}} \right)\frac{3{h^{\mbox{\tiny [1]}}}^{2}h^{\mbox{\tiny [2]}}}{32}+ \left( {\frac{5}{4}}+{\delta}+{\frac {\epsilon}{2}}+{\frac {3 {\delta}^{2}}{2}}+{\frac {3 \epsilon \delta}{2}}+{\frac {5 {\epsilon}^{2}}{16}}+ {5 {\delta}^{3}}+{\frac {15 \epsilon {\delta}^{2}}{2}}+{\frac {25 {\epsilon}^{2}\delta}{8}}+{\frac {19 {\epsilon}^{3}}{64}}-{\frac {35 {\delta}^{4}}{4}}-{\frac {35 \epsilon {\delta}^{3}}{2}}-{\frac {175 {\epsilon}^{2}{\delta}^{2}}{16}}-{\frac {133 {\epsilon}^{3}\delta}{64}}+{\frac {43 {\epsilon}^{4}}{1536}} \right) \frac{{h^{\mbox{\tiny [1]}}}^{4}}{32}=0.
\end{dmath*}
With $N=3,$ Galerkin projections based system of Eqns.(\ref{galerkin_lin_math}) is
\begin{dmath*}
\int_{0}^{4\pi}\cos(\tfrac{t}{2}) R_3(t) dt = 
  \frac{-1}{4}+\delta+\frac{\epsilon}{2}+ \left( -{\frac {1}{4}}+\delta \right) \frac{3 \epsilon h^{\mbox{\tiny [1]}}}{16}+ \left( -{\frac {1}{4}}+\delta \right) \frac{\epsilon h^{\mbox{\tiny [2]}}}{16}+ \left( -{\frac {1}{4}}+{\delta} \right) \frac{\epsilon h^{\mbox{\tiny [3]}}}{96}+ \left( -{\frac {3 \delta}{4}}-\frac {29 {\epsilon}}{64}+3 {\delta}^{2}+{\frac {29 {\epsilon} \delta}{16}}-\frac {{\epsilon}^{2}}{32} \right) \frac{\epsilon {h^{\mbox{\tiny [1]}}}^{2}}{16}+ \left( -{\frac {\delta}{4}}-{\frac {29 {\epsilon}}{192}}+{\delta}^{2}+{\frac {29 {\epsilon}\delta}{48}}-{\frac {{\epsilon}^{2}}{96}} \right) \frac{\epsilon h^{\mbox{\tiny [1]}} h^{\mbox{\tiny [2]}}}{16}+ \left( -{\frac {{\delta}^{2}}{4}}-{\frac {29 {\epsilon}\delta}{96}}-{\frac {1529 {\epsilon}^{2}}{18432}}+{\delta}^{3}+{\frac {29 {\epsilon} {\delta}^{2}}{24}}+{\frac {1433 {\epsilon}^{2}\delta}{4608}}-{\frac {109 {\epsilon}^{3}}{9216}} \right) \frac{\epsilon {h^{\mbox{\tiny [1]}}}^{3}}{16}=0,
\end{dmath*} 

\begin{dmath*}
\int_{0}^{4\pi}\cos(\tfrac{3t}{2}) R_3(t) dt = 
   \frac{\epsilon}{2}+\left( {\frac {27}{4}}-{3 \delta}+{\frac {3 {\epsilon}}{2}} \right) \frac{\epsilon h^{\mbox{\tiny [1]}}}{16}+ \left( {\frac {9}{4}}-\delta+\frac{\epsilon}{2} \right)\frac{\epsilon  h^{\mbox{\tiny [2]}}}{16}+ \left( {\frac {3}{8}}-{\frac {\delta}{6}}+{\frac {{\epsilon}}{12}} \right) \frac{\epsilon h^{\mbox{\tiny [3]}}}{16}+ \left( {\frac {27 \delta}{4}}+{\frac {135 {\epsilon}}{32}}-3 {\delta}^{2}-{\frac {3 {\epsilon} \delta}{8}}+{\frac {15 {\epsilon}^{2}}{16}} \right) \frac{\epsilon {h^{\mbox{\tiny [1]}}}^{2}}{16}+ \left( {\frac {9 \delta}{4}}+{\frac {45 {\epsilon}}{32}}-{\delta}^{2}-{\frac {{\epsilon} \delta}{8}}+{\frac {5 {\epsilon}^{2}}{16}} \right) \frac{\epsilon h^{\mbox{\tiny [1]}} h^{\mbox{\tiny [2]}}}{16}+ \left( {\frac {9 {\delta}^{2}}{4}}+{\frac {45 {\epsilon} \delta}{16}}+{\frac {819 {\epsilon}^{2}}{1024}}-{\delta}^{3}-{\frac {3 {\epsilon}{\delta}^{2}}{4}}+{\frac {69 {\epsilon}^{2}\delta}{256}}+{\frac {1639 {\epsilon}^{3}}{9216}} \right)\frac{ \epsilon {h^{\mbox{\tiny [1]}}}^{3}}{16}=0,
\end{dmath*} 

\begin{dmath*}
\hspace*{-5mm} \int_{0}^{4\pi}\cos(\tfrac{5t}{2}) R_3(t) dt = 
   -{\frac {3 {\epsilon}^{2} h^{\mbox{\tiny [1]}}}{32}}-\frac{ {\epsilon}^{2} h^{\mbox{\tiny [2]}}}{32}-{\frac {{\epsilon}^{2} h^{\mbox{\tiny [3]}}}{192}}+\left( -{\frac {25}{4}}-{{23 \delta}}-{15 {\epsilon}} \right) \frac{\epsilon {h^{\mbox{\tiny [1]}}}^{2}}{256}+ \left( -{\frac {25}{12}}-{\frac {23 \delta}{3}}-{{5 {\epsilon}}} \right) \frac{{\epsilon}^{2} h^{\mbox{\tiny [1]}} h^{\mbox{\tiny [2]}}}{256}+ \left( -{\frac {25 \delta}{4}}-{\frac {1375 {\epsilon}}{384}}-{ {11 {\delta}^{2}}}-{\frac {1385 {\epsilon} \delta}{96}}-{\frac {1639 {\epsilon}^{2}}{384}} \right) \frac{{\epsilon}^{2} {h^{\mbox{\tiny [1]}}}^{3}}{384}=0
\end{dmath*} 
and
\begin{flalign*}
\int_{0}^{4\pi}\cos(\tfrac{7t}{2}) R_3(t) dt = 
   \frac{{\epsilon}^{3} {h^{\mbox{\tiny [1]}}}^{2}}{512}+\frac{{\epsilon}^{3}{h^{\mbox{\tiny [1]}}}{h^{\mbox{\tiny [2]}}}}{1536}+\left(\frac{49}{4}+95 {\delta}+{55 \epsilon} \right) \frac{{\epsilon}^{3} {h^{\mbox{\tiny [1]}}}^{3}}{73728}=0.
\end{flalign*}

\subsection{Mathieu equation with impulsive parametric excitation}
\label{lin_impulse_eqns}
	
Considering Eqn.(\ref{taylor_approx1_impulse}) with $N=3$ and substituting expressions for $\lambda^{\tiny [1]},\ \lambda^{\tiny [2]},\ \lambda^{\tiny [3]}$ and $\lambda^{\tiny [4]}$ gives
\begin{dmath*}
  \left( -1+ \Delta+ {\frac {\epsilon}{\pi}} \right)2 h^{\mbox{\tiny [1]}}+ \left( -1+\Delta+{\frac {\epsilon}{\pi}} \right)\frac{3 h^{\mbox{\tiny [2]}}}{4}+ \left( -1+\Delta+{\frac {\epsilon}{\pi}} \right) \frac{h^{\mbox{\tiny [3]}}}{6}+ \left( -1+\Delta+{\frac {\epsilon}{\pi}} \right) \frac{h^{\mbox{\tiny [4]}}}{48}+ \left(- \frac{3}{2}+\Delta+{\frac {\epsilon}{\pi}}+\frac{{\Delta}^{2}}{2}+{\frac {\Delta \epsilon}{\pi}} \right)\frac{3  {h^{\mbox{\tiny [1]}}}^{2}}{2}+ \left( -{\frac{3}{4}}+\Delta+{\frac {\epsilon}{\pi}}+\frac{{\Delta}^{2}}{2}+\frac{\Delta \epsilon}{2\pi} \right)\frac{3  h^{\mbox{\tiny [1]}}h^{\mbox{\tiny [2]}}}{2}+ \left( -{\frac{3}{2}}+\Delta+{\frac {\epsilon}{\pi}}+\frac{{\Delta}^{2}}{2}+{\frac {\Delta \epsilon}{\pi}} \right)\frac{{h^{\mbox{\tiny [2]}}}^{2}}{16}+ \left( -\frac{1}{2}+\frac{\Delta}{3}+{\frac {\epsilon}{3\pi}}+\frac{{\Delta}^{2}}{6}+{\frac {\Delta \epsilon}{3\pi}} \right)\frac{ h^{\mbox{\tiny [1]}}h^{\mbox{\tiny [3]}}}{4}+ \left( -5+3\Delta+{\frac {3\epsilon}{\pi}}+{\Delta}^{2}+{\frac {2 \Delta \epsilon}{\pi}}+{\Delta}^{3}+{\frac {{3\Delta}^{2}\epsilon}{\pi}}-{\frac {{2\epsilon}^{3}}{3\pi}} \right) \frac{{h^{\mbox{\tiny [1]}}}^{3}}{4}+ \left( -{\frac{15}{2}}+{\frac {9 \Delta}{2}}+{\frac {9 \epsilon}{2 \pi}}+{\frac {3 {\Delta}^{2}}{2}}+{\frac {3\Delta \epsilon}{\pi}}+{\frac {3 {\Delta}^{3}}{2}}+{\frac {9 {\Delta}^{2}\epsilon}{2 \pi}}-{\frac {{\epsilon}^{3}}{\pi}} \right) \frac{{h^{\mbox{\tiny [1]}}}^{2} h^{\mbox{\tiny [2]}}}{16}+ \left( -{\frac{35}{4}}+{ {5 \Delta}}+{\frac {5 \epsilon}{\pi}}+{\frac {3 {\Delta}^{2}}{2}}+{\frac {3 \Delta \epsilon}{\pi}}+{\Delta}^{3}+{\frac {3 {\Delta}^{2}\epsilon}{\pi}}-{\frac {{2\epsilon}^{3}}{3\pi}}+{\frac {5 {\Delta}^{4}}{4}}+{\frac {5 {\Delta}^{3}\epsilon}{\pi}}-{\frac {10 \Delta {\epsilon}^{3}}{3 \pi}}-{\frac {{4\epsilon}^{4}}{{3\pi}^{2}}} \right) \frac{{h^{\mbox{\tiny [1]}}}^{4}}{32}=0.
\end{dmath*} 
With $N=3,$ Galerkin projections based system of Eqns. is
\begin{dmath*} 
\int_{0}^{2\pi}\cos \left( t \right) R_3(t) dt = 
 -\pi+\pi \Delta+\epsilon+\left( -1+\Delta \right)\frac{3 \epsilon  h^{\mbox{\tiny [1]}}}{4}+ \left( -1+\Delta \right) \frac{\epsilon h^{\mbox{\tiny [2]}}}{4}+ \left( -1+\Delta \right) \frac{\epsilon h^{\mbox{\tiny [3]}}}{24}+\left( -\Delta+{\frac { \left( 2 {\pi}^{2}-3 \right){\epsilon}^{2}}{12 {\pi}}}+{\Delta}^{2}+{\frac { \left( -2 {\pi}^{2}+3\right) {\epsilon}\Delta}{12 {\pi}}}-\frac{{\epsilon}^{2}}{2} \right) \frac{3 {h^{\mbox{\tiny [1]}}}^{2}}{4}+ \left( -{ \epsilon \Delta}+{\frac { \left( 2 {\pi}^{2}-3 \right) {\epsilon}^{2}}{12 {\pi}}}+\epsilon {\Delta}^{2}+{\frac { \left( -2 {\pi}^{2}+3 \right) {\epsilon}^{2}\Delta}{12 {\pi}}}-\frac{{\epsilon}^{3}}{2} \right) \frac{h^{\mbox{\tiny [1]}}h^{\mbox{\tiny [2]}}}{4}+ \left( -{\Delta}^{2}+{\frac { \left( 2 {\pi}^{2}-3\right) {\epsilon}\Delta}{6 {\pi}}}+{\frac { \left( 14 {\pi}^{2}-3 \right) {\epsilon}^{2}}{48 {\pi}^{2}}}+{\Delta}^{3}+{\frac { \left( -2 {\pi}^{2}+3\right) {\epsilon}^{2}{\Delta}^{2}}{6 {\pi}}}+{\frac { \left( -62 {\pi}^{2}+3 \right) {\epsilon}^{2}\Delta}{48 {\pi}^{2}}}-{\frac {{\epsilon}^{4}}{2\pi}} \right)\frac{ {h^{\mbox{\tiny [1]}}}^{3}}{4}=0,
\end{dmath*} 

\begin{dmath*}
\int_{0}^{2\pi}{\it H} \left(t -\pi\right) \sin \left( t \right) R_3(t) dt = 
  \left( -1+\Delta \right) \frac{3\pi \epsilon h^{\mbox{\tiny [1]}}}{8}+ \left( -1+\Delta \right) \frac{\pi \epsilon h^{\mbox{\tiny [2]}}}{8}+ \left( -1+\Delta \right)\frac{ \pi \epsilon h^{\mbox{\tiny [3]}}}{48}+ \left( -\pi \Delta-\frac {5 {\epsilon}}{4}+\pi {\Delta}^{2}+\frac{\epsilon\Delta}{4} \right) \frac{3 \epsilon {h^{\mbox{\tiny [1]}}}^{2}}{8}+ \left( -\pi \Delta-\frac {5 {\epsilon}}{4}+\pi {\Delta}^{2}+\frac{\epsilon\Delta}{4} \right) \frac{\epsilon h^{\mbox{\tiny [1]}}h^{\mbox{\tiny [2]}}}{8}+ \left( \pi {\Delta}^{2}-{\frac {5 {\epsilon}\Delta}{2}}+{\frac { \left( {{5 {\pi}^{2}}}-{{51}} \right) {\epsilon}^{2}}{6\pi}}+\pi {\Delta}^{3}+8{\epsilon}{\Delta}^{2}+\frac { \left( -{ {5 {\pi}^{2}}}+3 \right) {\epsilon}^{2}\Delta}{48\pi} \right) \frac{\epsilon{h^{\mbox{\tiny [1]}}}^{3}}{8}=0,
\end{dmath*} 

\begin{dmath*}
\int_{0}^{2\pi}{\it H} \left(t -\pi\right) \cos \left( t \right) R_3(t) dt = 
 -\frac{\pi}{2}+\epsilon+\frac{\pi \Delta}{2}+ \left( {3}+\Delta \right) \frac{3 \epsilon h^{\mbox{\tiny [1]}}}{8}+ \left( 3+\Delta \right) \frac{\epsilon h^{\mbox{\tiny [2]}}}{8}+ \left( 1+\frac{\Delta}{3} \right) \frac{\epsilon h^{\mbox{\tiny [3]}}}{16}+ \left( { {9 \Delta}}+{\frac { \left( 2 {\pi}^{2}+45 \right) {\epsilon}}{4 {\pi}}}+3 {\Delta}^{2}+{\frac { \left( -2 {\pi}^{2}+3  \right) {\epsilon}\Delta}{4 {\pi}}}-{3 {\epsilon}^{2}} \right) \frac{\epsilon {h^{\mbox{\tiny [1]}}}^{2}}{8}+ \left( {3 \Delta}+{\frac { \left( {\pi}^{2}+25 \right) {\epsilon}}{6 {\pi}}}+{{\Delta}^{2}}+{\frac { \left( -2 {\pi}^{2}+3 \right) {\epsilon}\Delta}{12 {\pi}}}-{{\epsilon}^{2}} \right) \frac{\epsilon h^{\mbox{\tiny [1]}}h^{\mbox{\tiny [2]}}}{8}+ \left( {3 {\Delta}^{2}}+{\frac { \left( 2 {\pi}^{2}+45  \right) \epsilon\Delta}{6 {\pi}}}+{\frac { \left( -10 {\pi}^{2}+141 \right) {\epsilon}^{2}}{48 {\pi}^{2}}}+{{\Delta}^{3}}+{\frac { \left( -2 {\pi}^{2}+3 \right) {\epsilon}{\Delta}^{2}}{6 {\pi}}}+{\frac { \left( -110 {\pi}^{2}+3 \right) {\epsilon}^{2}\Delta}{48 {\pi}^{2}}}-{\frac {{\epsilon}^{3}}{\pi}} \right) \frac{\epsilon {h^{\mbox{\tiny [1]}}}^{3}}{8}=0\quad 
\end{dmath*} 
and
\begin{dmath*}
\int_{0}^{2\pi}tR_3(t) dt = 
 -\pi \epsilon-\left(3 h^{\mbox{\tiny [1]}}-h^{\mbox{\tiny [2]}}-\frac{h^{\mbox{\tiny [3]}}}{6}\right)\pi \epsilon \Delta+ \left( -{\Delta}^{2}+\frac{{\epsilon}^{2}}{8} \right) 3\pi\epsilon{h^{\mbox{\tiny [1]}}}^{2}+ \left( -{\Delta}^{2}+\frac{{\epsilon}^{2}}{8} \right)\pi\epsilon h^{\mbox{\tiny [1]}} h^{\mbox{\tiny [2]}}+ \left( -{\Delta}^{3}+\frac{3 {\epsilon}^{2}\Delta}{8}+\frac{{\epsilon}^{3}}{\pi} \right)\pi\epsilon {h^{\mbox{\tiny [1]}}}^{3}=0.
\end{dmath*}
\clearpage

\section{$4\pi$-periodic solution via HAM for Mathieu equation with impulsive parametric excitation}
\label{lin_impulse_sol}
	
The solution to Eqn.(\ref{ode_lin_impulse_math}) $\hat{x}_3(t),$ periodic with period $4\pi,$ approximated using HAM is

\begin{dmath}\label{x_sol_2}
\hat{x}_3(t) =
-\epsilon \left( {\frac {\left( -48\,{\delta}^{2}+{\epsilon}^{2} \right) {h^{\mbox{\tiny [1]}}}^{3}}{48\,{\pi }}}+{\frac {\delta{h^{\mbox{\tiny [1]}}}^{2}}{\pi }}+{\frac {\delta h^{\mbox{\tiny [1]}}h^{\mbox{\tiny [2]}}}{\pi }}+{\frac {3h^{\mbox{\tiny [1]}}}{\pi }}+{\frac{h^{\mbox{\tiny [2]}}}{\pi }}+{\frac {h^{\mbox{\tiny [3]}}}{6\pi }}-2 \right)\sin(\tfrac{t}{2})
-\frac{{\epsilon}^{2}h^{\mbox{\tiny [1]}}}{4}\left( {\frac {\left( 4{\pi }\Delta+{\epsilon} \right) {h^{\mbox{\tiny [1]}}}^{2}}{2}}+{3{h^{\mbox{\tiny [1]}}}}+h^{\mbox{\tiny [2]}} \right) \sin(\tfrac{t}{2}) {\it H} \left( t-\pi \right)
 -\frac{3{\epsilon}^{2}h^{\mbox{\tiny [1]}}}{4}\left( {\frac {\left(4{\pi }\Delta+{\epsilon} \right) {h^{\mbox{\tiny [1]}}}^{2}}{2{\pi }}}+{3{h^{\mbox{\tiny [1]}}}}+{h^{\mbox{\tiny [2]}}} \right) \sin(\tfrac{t}{2}) H \left( t-3\pi \right)
-\left( {\frac {\epsilon \left( 16\,{\Delta}^{2}-{\epsilon}^{2} \right) {h^{\mbox{\tiny [1]}}}^{3}}{16}}+3{h^{\mbox{\tiny [1]}}}
^{2}\epsilon\,\Delta+\epsilon\,\Delta\,h^{\mbox{\tiny [1]}}h^{\mbox{\tiny [2]}}+3\,\epsilon\,h^{\mbox{\tiny [1]}}+\epsilon\,h^{\mbox{\tiny [2]}}+\frac{\epsilon\,h^{\mbox{\tiny [3]}} }{6}\right) \cos \left( \tfrac{t}{2}\right) {\it H} \left( t-\pi \right)
 -  \epsilon\left( {\frac { \left( 16{\Delta}^{2}-5{\epsilon}^{2} \right) {h^{\mbox{\tiny [1]}}}^{3}}{16}}+3\Delta{h^{\mbox{\tiny [1]}}}^{2}+\Delta\,h^{\mbox{\tiny [1]}}h^{\mbox{\tiny [2]}}+ 3h^{\mbox{\tiny [1]}}+h^{\mbox{\tiny [2]}}+\frac{h^{\mbox{\tiny [3]}}}{6} \right) \cos(\tfrac{t}{2}) {\it H} \left( t-3\pi \right)  
 +\frac{\epsilon}{2\pi}\left( {\frac { \left( 48{\Delta}^{2}-{
\epsilon} \right) {h^{\mbox{\tiny [1]}}}^{3}}{48}}+
{3{h^{\mbox{\tiny [1]}}}^{2}\Delta}+{\Delta\,h^{\mbox{\tiny [1]}}h^{\mbox{\tiny [2]}}}+{3h^{\mbox{\tiny [1]}}}+{h^{\mbox{\tiny [2]}}}+{\frac {h^{\mbox{\tiny [3]}}}{6}}
 \right)t \cos(\tfrac{t}{2}) 
+\frac{{\epsilon}^{2}h^{\mbox{\tiny [1]}}}{4\pi}\left( {\frac { \left( 4 {\pi }\Delta+{\epsilon} \right) {h^{\mbox{\tiny [1]}}}^{2}}{2\pi }}+{{3{h^{\mbox{\tiny [1]}}}}}+{h^{\mbox{\tiny [2]}}} \right)t \sin(\tfrac{t}{2}) {\it H}\left( t-\pi \right) 
+\frac{{\epsilon}^{2}h^{\mbox{\tiny [1]}}}{4\pi}\left( {\frac { \left( 4 {\pi }\Delta+ 
{\epsilon} \right) {h^{\mbox{\tiny [1]}}}^{2}}{2\pi }}+{{3
{h^{\mbox{\tiny [1]}}}}}+{{h^{\mbox{\tiny [2]}}}} \right) t\sin(\tfrac{t}{2}) {\it H}\left( t-3\pi \right)
-{\frac {{\epsilon}^{3}{h^{\mbox{\tiny [1]}}}^{3}}{16\pi }}t\cos(\tfrac{t}{2}) {\it H} \left( t-\pi \right)
-{\frac {3{\epsilon}^{3}{h^{\mbox{\tiny [1]}}}^{3} }{16\pi }}t \cos(\tfrac{t}{2}) {\it H} \left( t-3\pi \right)
- \frac{{\epsilon}^{2}h^{\mbox{\tiny [1]}}}{16\pi^2}\left( {\frac {\left(4\,\Delta+{\epsilon} \right) {h^{\mbox{\tiny [1]}}}^{2}}{2}}+{h^{\mbox{\tiny [1]}}}+h^{\mbox{\tiny [2]}} \right){t}^{2} \sin(\tfrac{t}{2})
+{\frac {{\epsilon}^{3}{h^{\mbox{\tiny [1]}}}^{3}}{32{\pi }^{2}}}{t}^{2} \cos(\tfrac{t}{2}) {\it H} \left( t-\pi \right)
+{\frac {{\epsilon}^{3}{h^{\mbox{\tiny [1]}}}^{3} }{32{\pi }^{2}}}{t}^{2}\cos(\tfrac{t}{2}) {\it H} \left( t-3\pi\right) 
-{\frac {{\epsilon }^{3}{h^{\mbox{\tiny [1]}}}^{3}}{192\,{\pi }^{3}}}{t}^{3}\cos(\tfrac{t}{2}) 
. 
\end{dmath}

\end{appendices}

%\printbibliography[omitnumbers=false]
%\bibliographystyle{unsrt}

\end{document}